 \loadeufm \loadeusm \loadbold
\documentstyle{amsppt}
\NoRunningHeads \magnification=1200 \pagewidth{32pc}
\pageheight{42pc} \vcorrection{1.2pc}

\define\wh{\widehat}

\define\bc{\Bbb C}
\define\bz{\Bbb Z}
\define\br{\Bbb R}
\define\bn{\Bbb N}

\define\bda{\boldkey A}
\define\bdb{\boldkey B}

\define\bbda{\bar\boldkey A}
\define\bbdb{\bar\boldkey B}
\define\bk{\Bbb K}
\define\eu{e^{(\mu)}}
\document

\topmatter
\title Unitary representations of the extended affine Lie algebra
$\widetilde{\frak{gl}_{3}(\bc_q)}$\endtitle
\author Ziting Zeng
\endauthor
\address
Department of Mathematics and Statistics, York University,
Toronto, \newline Canada M3J 1P3
\endaddress
\email \ ziting\@mathstat.yorku.ca
\endemail

\abstract We present modules for the extended affine Lie algebra
$\widetilde{\frak{gl}_{3}(\bc_q)}$ by using the idea of free
fields. A necessary and sufficient condition for the modules being
unitary is given.
\endabstract

\endtopmatter
\document

\subhead \S 0. Introduction \endsubhead

Let $q$ be a non-zero complex number. A quantum $2$-torus
associated to $q$ (see [M]) is the unital associative
$\bc$-algebra $\bc_{q}[s^{\pm 1}, t^{\pm 1}]$ (or, simply
$\bc_{q}$) with generators $s^{\pm 1},  t^{\pm 1}$ and relations
$$s s^{-1} =s^{-1}s =t t^{-1}= t^{-1}t=1 \, \text{  and } \, \ ts = q st. \tag 0.1$$

Define $\kappa, : \bc_q \to \bc$ to be  a $\bc$-linear function
given by
$$\kappa(s^m t^n) =\delta_{(m,n),(0,0)} \tag 0.2$$

Let $d_s$, $d_t$ be the degree operators on $\bc_q$ defined by
$$d_s(s^mt^n) = ms^m t^n, \, \,  d_t(s^m t^n) = n s^m t^n \tag 0.3$$
 for $m, n\in \bz$.

\medskip

Let $\frak{gl}_{3}(\bc_q)$ be the Lie algebra of $3$ by $3$
matrices whose entries are from $\bc_q$. We form a  natural
central extension of $\frak{gl}_{3}(\bc_q)$ as follows.
$$\widehat{\frak{gl}_{3}(\bc_q)} = \frak{gl}_{3}(\bc_q)\oplus \bc c_s\oplus
\bc c_t \tag 0.4$$ with Lie bracket
 $$\align & [E_{ij}(s^{m_1}t^{n_1}), E_{kl}(s^{m_2}t^{n_2})]\tag 0.5\\
=&\delta_{jk}q^{n_1m_2}E_{il}(s^{m_1+m_2}t^{n_1+n_2})-\delta_{il}q^{n_2
m_1}E_{kj}
(s^{m_1+m_2}t^{n_1+n_2})\\
&+m_1q^{n_1 m_2}\delta_{jk}\delta_{il}\delta_{m_1+m_2,
0}\delta_{n_1 + n_2, 0}c_s
+n_1q^{n_1m_2}\delta_{jk}\delta_{il}\delta_{m_1+m_2,
0}\delta_{n_1+n_2, 0}c_t
\endalign$$
for $m_1, m_2, n_1, n_2\in \bz$, $1\leq i, j, k, l\leq 3$, where
$E_{ij}$ is the matrix whose $(i, j)$-entry is $1$ and $0$
elsewhere, and $c_s$ and $c_t$ are central elements of
$\widehat{\frak{gl}_{3}(\bc_q)}$.

The derivations $d_s$ and $d_t$ can be extended to derivations on
$\frak{gl}_{3}(\bc_q)$. Now we can define the semi-direct product
of the Lie algebra $\widehat{\frak{gl}_{3}(\bc_q)}$ and those
derivations:
$$\widetilde{\frak{gl}_{3}(\bc_q)} =\widehat{\frak{gl}_{3}(\bc_q)} \oplus \bc d_s
\oplus \bc d_t. \tag 0.6$$

The Lie algebra $\widetilde{\frak{gl}_{3}(\bc_q)}$ is  an extended
affine
 Lie algebra of type $A_{2}$ with nullity $2$. (See [AABGP] and
[BGK] for definitions).

\medskip

Extended affine Lie algebras are a higher dimensional
generalization of affine Kac-Moody Lie algebras introduced by
[H-KT] and systematically studied in [AABGP] and [BGK]. It turns
out that any extended affine Lie algebra of type $A$ is
coordinated by a quantum torus (or a nonassociative torus for some
small rank cases).  Representations for extended affine Lie
algebras coordinated by quantum tori and Lie algebras  related to
quantum tori have been studied in [JK2], [BS], [G1,2,3], [ER1,2],
[EB],[GZ] [EZ], [LT1,2], [G-KK], [VV], [Mi], [ZZ], [BZ], [SZ], [L]
and [BEG], and among others.

\medskip

The Wakimoto's  free fields construction provides a remarkable way
to realize  affine Kac-Moody Lie algebras (see [W2], [FF] and
[EFK]). In [GZ], we used Wakimoto's idea to construct a class of
representations for $\widetilde{\frak{gl}_{2}(\bc_q)}$ and found
out the necessary and sufficient condition for the representations
being unitary. In this paper, we will continue to construct
representations for $\widetilde{\frak{gl}_{3}(\bc_q)}$.  As
witnessed in [FF], the  realization for
$\widetilde{\frak{gl}_{3}(\bc_q)}$ is much more subtle and
complicated than the one for $\widetilde{\frak{gl}_{2}(\bc_q)}$.
We then go on to construct a hermitian form and to determine when
the form is positive definite (so the representations are
unitary). Unlike [GZ] in which we defined the form on the monomial
basis for the module (this idea goes back to [W1]), we define the
form directly on the basis consisting of certain iterated module
actions on a ``highest weight vector" $1$. This facilitates the
verification of the defined form  being a hermitian from.

\medskip

{\it Throughout this paper, we denote the field of complex
numbers, real numbers and
 the ring of integers by $\bc$, $\br$ and $\bz$ respectively.}

\medskip

 \subhead \S 1. Module for
$\widetilde{\frak{gl}_{3}(\bc_q)}$\endsubhead

In this section, we use Wakimoto's idea [W1] to
 construct
$\widetilde{\frak{gl}_{3}(\bc_q)}$-modules as was done in [GZ].

 Let
$\bk_1=\{(3m+1,3n+1),m,n\in\bz\}$, and
$\bk_{-1}=\{(3m-1,3n-1),m,n\in\bz\}$. If $\bda=(3m+1,3n+1)\in
\bk_1$, we always write $\bda_1=m, \bda_2=n$, and similarly, if
$\bdb=(3m-1,3n-1)\in \bk_{-1}$, then $\bdb_1=m, \bdb_2=n$. Set
$$V = \bc[x_{\bda},x_{\bdb}: \bda\in\bk_1, \bdb\in\bk_{-1}]\tag 1.1$$
 be the (commutative) polynomial ring of infinitely many variables.
The operators $x_{(m, n)}$ and $\frac{\partial}{\partial
x_{(m,n)}}$ act on $V$ as the usual multiplication and
differentiation operators respectively.

\medskip

Given a family of $2\times 2$ lower triangular matrices
$$X_{m,n}=\pmatrix a_{(m, n)} & 0 \\
c_{(m, n)} & d_{(m, n)}\endpmatrix \in \text{ SL}_2(\bc) $$ for
$(m,n)\in \bk_1\bigcup\bk_{-1}$ (so $a_{(m,n)} d_{(m, n)}=1$).

 Set
 $$\align & P_{(m,n)}=a_{(m,n)}\frac{\partial}{\partial
 x_{(m,n)}}\tag 1.2\\
& Q_{(m,n)}=c_{(m,n)}\frac{\partial}{\partial
x_{(m,n)}}+d_{(m,n)}x_{(m,n)}\tag 1.3\endalign $$ for $(m, n)\in
\bk_1\bigcup\bk_{-1}$. Then for $\bda, \bda^\prime \in \bk_1,
\bdb,\bdb^\prime \in \bk_{-1}$,
$$\align & [P_\bda,P_{\bda^\prime}]=[Q_\bda, Q_{\bda^\prime}]=[P_\bdb, P_{\bdb^\prime}]=[Q_\bdb, Q_{\bdb^\prime}]=0\\
&[P_\bda, P_\bdb] = [P_\bda, Q_\bdb]=[Q_\bda, Q_\bdb]=[P_\bdb,
Q_\bda]=0\\
&[P_\bda, Q_{\bda^\prime}] = \delta_{\bda, \bda^\prime}, \ \
[P_\bdb, Q_{\bdb^\prime}] = \delta_{\bdb, \bdb^\prime}.
 \endalign $$

Fix a complex number $\mu$, we  define the following operators on
 $V$:
$$\align  e_{21}^{(\mu)}(m_1, n_1) =& -q^{-m_1n_1}\mu P_{-(3m_1-1,3n_1-1)}\tag 1.4\\
&-\sum\Sb\bda,\bda'\in\bk_1 \endSb q^{n_1\bda_1'+ \bda_2 m_1
+\bda_2\bda_1'} Q_{\bda+\bda'+(3m_1-1, 3n_1-1)}
P_{\bda}P_{\bda'}\\
&-\sum\Sb\bda\in\bk_1\\\bdb\in\bk_{-1} \endSb q^{n_1\bda_1+ \bdb_2
m_1 +\bdb_2\bda_1} Q_{\bda+\bdb+(3m_1-1, 3n_1-1)}
P_{\bda}P_{\bdb}\\
e_{12}^{(\mu)}(m_1, n_1)=& Q_{(3m_1+1,3n_1+1)}\tag 1.5\\
e_{11}^{(\mu)}(m_1, n_1)=&\sum_{\bda\in\bk_1}
q^{\bda_1n_1}Q_{(3m_1, 3n_1)+\bda}
P_{\bda}+\frac{1}{2}\mu\delta_{(m_1,n_1),(0,0)}\tag1.6\\
e_{22}^{(\mu)}(m_1, n_1)=& -\sum_{\bda\in\bk_1} q^{\bda_2m_1}
Q_{(3m_1, 3n_1)+\bda} P_{\bda}\tag 1.7\\
&-\sum_{\bdb\in\bk_{-1}} q^{\bdb_2m_1} Q_{(3m_1, 3n_1)+\bdb}
P_{\bdb}-\frac{1}{2}\mu\delta_{(m_1,n_1),(0,0)}\endalign$$

$$\align
e_{23}^{(\mu)}(m_1, n_1) &= -q^{-m_1n_1}\mu P_{-(3m_1+1,3n_1+1)}\tag 1.8\\
&-\sum\Sb\bda\in\bk_1\\\bdb\in\bk_{-1} \endSb q^{n_1\bdb_1+ \bda_2
m_1 +\bda_2\bdb_1} Q_{\bda+\bdb+(3m_1+1, 3n_1+1)}
P_{\bda}P_{\bdb}\\
&-\sum\Sb\bdb,\bdb'\in\bk_{-1} \endSb q^{n_1\bdb_1'+ \bdb_2 m_1
+\bdb_2\bdb_1'} Q_{\bdb+\bdb'+(3m_1+1, 3n_1+1)}
P_{\bdb}P_{\bdb'}\\
e_{32}^{(\mu)}(m_1, n_1)&= Q_{(3m_1-1,3n_1-1)}\tag 1.9\\
e_{31}^{(\mu)}(m_1, n_1)&=\sum_{\bda\in\bk_1}
q^{\bda_1n_1}Q_{(3m_1-2, 3n_1-2)+\bda}
P_{\bda}\tag 1.10\\
e_{13}^{(\mu)}(m_1, n_1)&=\sum_{\bdb\in\bk_{-1}}
q^{\bdb_1n_1}Q_{(3m_1+2, 3n_1+2)+\bdb}
P_{\bdb}\tag 1.11\\
e_{33}^{(\mu)}(m_1, n_1)&=\sum_{\bdb\in\bk_{-1}}
q^{\bdb_1n_1}Q_{(3m_1, 3n_1)+\bdb}
P_{\bdb}+\frac{1}{2}\mu\delta_{(m_1,n_1),(0,0)}\tag 1.12\\
D_1 ^{(\mu)}&= \sum_{\bda\in\bk_1}\bda_1 Q_{\bda}P_{\bda}+\sum_{\bdb\in\bk_{-1}}\bdb_1 Q_{\bdb}P_{\bdb}\tag 1.13\\
D_2 ^{(\mu)}&= \sum_{\bda\in\bk_1}\bda_2
Q_{\bda}P_{\bda}+\sum_{\bdb\in\bk_{-1}}\bdb_2 Q_{\bdb}P_{\bdb}\tag
1.14
\endalign$$

Although the operators are infinite sums, they are well-defined as
operators on $V$. Now we have the following result:
\proclaim{Theorem 1.15} The linear map $\pi:
\widetilde{\frak{gl}_{3}(\bc_q)} \to \text{End }V$ given by
$$\align & \pi(E_{ij}(s^{m_1}t^{n_1})) = \eu_{ij}(m_1, n_1),\\
&\pi(d_s) = D_1^{(\mu)}, \quad \pi(d_t) = D_2^{(\mu)},\quad
\pi(c_s) = \pi(c_t) = 0,
\endalign $$
for $m_1, n_1\in \bz, 1\leq i, j\leq 3$, is a Lie algebra
homomorphism.
\endproclaim
\demo{Proof} The proof is straightforward. However, we would like
to provide a few details. We shall do this orderly so that we
won't miss any case.

First, we have
$$\align & [\eu_{11}(m_1, n_1), \eu_{11}(m_2, n_2)]\\
=&[\sum_{\bda\in\bk_1} q^{\bda_1n_1}Q_{(3m_1, 3n_1)+\bda}
P_{\bda}, \sum_{\bda'\in\bk_1} q^{\bda'_1n_1}Q_{(3m_2,
3n_2)+\bda'} P_{\bda'}]\\
=&\sum_{\bda'\in\bk_1}q^{(m_2+\bda'_1)n_1+\bda'n_2}Q_{(3m_1,3n_1)+(3m_2+3n_2)+\bda'}P_{\bda'}
+\frac{1}{2}q^{m_2n_1}\mu\delta_{(m_1+m_2,n_1+n_2),(0,0)}\\
&-\sum_{\bda\in\bk_1}q^{\bda_1n_1+(m_1+\bda_1)n_2}q_{(3m_1,3n_1)+(3m_2,3n_2)+\bda}P_{\bda}
-\frac{1}{2}q^{m_1n_2}\mu\delta_{(m_1+m_2,n_1+n_2),(0,0)}\\
=&q^{m_2n_1}\eu_{11}(m_1+m_2,n_1+n_2)-q^{m_1n_2}\eu_{11}(m_1+m_2,n_1+n_2).
\endalign$$
The following two brackets are easy.
$$ \align
&[\eu_{11}(m_1,n_1),\eu_{12}(m_2,n_2)]=q^{m_2n_1}\eu_{12}(m_1+m_2,n_1+n_2)\\
&[\eu_{11}(m_1,n_1),\eu_{13}(m_2,n_2)]=q^{m_2n_1}\eu_{13}(m_1+m_2,n_1+n_2).
\endalign $$
$$\align  & [\eu_{11}(m_1,n_1),\eu_{21}(m_2,n_2)]\\
=&\sum_{\bda\in\bk_1}(-\mu
)q^{\bda_1n_1-m_2n_2}[Q_{(3m_1,3n_1)+\bda}P_\bda,
P_{-(3m_2-1,3n_2-1)}]\\
&-\sum\Sb \bda\in\bk_1\\ \bbda, \bbda'\in\bk_1\endSb
q^{\bda_1n_1+\bbda'_1n_2+\bbda_2m_2+\bbda_2\bbda'_1}[Q_{(3m_1,3n_1)+\bda}P_\bda,
Q_{\bbda+\bbda'+(3m_2-1,3n_2-1)}P_{\bbda}P_{\bbda'}]\\
&-\sum\Sb \bda\in\bk_1\\ \bbda\in\bk_1\\ \bdb\in \bk_{-1}\endSb
q^{\bda_1n_1+n_2\bbda_1+\bdb_2m_2+\bdb_2\bbda_1}[Q_{(3m_1,3n_1)+\bda}P_\bda,
Q_{\bbda+\bdb+(3m_2-1,3n_2-1)}P_{\bbda}P_{\bdb}]\\
=&\mu
q^{-(m_1+m_2)(n_1+n_2)+m_1n_2}P_{-(3(m_1+m_2)-1,3(n_1+n_2)-1)}\\
&-\sum_{\bbda,\bbda'\in\bk_1}q^{(\bbda_1+\bbda'_1+m_2)n_1+n_2\bbda'_1+\bbda_2m_2+\bbda_2\bbda'_1}
Q_{(3m_1,3n_1)+\bbda+\bbda'+(3m_2-1,3n_2-1)}P_{\bbda} P_{\bbda'}\\
&+\sum_{\bda,\bbda\in\bk_1}q^{\bda_1n_1+n_2(m_1+\bda_1)+\bbda_2m_2+\bbda_2(m_1+\bda_1)}
Q_{\bbda+(3m_1,3n_1)+\bda+(3m_2-1,3n_2-1)}P_{\bbda}P_{\bda}\\
&+\sum_{\bda,\bbda'\in\bk_1}q^{\bda_1n_1+n_2\bda'_1+(n_1+\bda_2)m_2+(n_1+\bda_2)\bbda'_1}
Q_{\bda+(3m_1,3n_1)+\bbda'+(3m_2-1,3n_2-1)}P_{\bbda'}P_{\bda}\\
&+\sum\Sb \bda\in\bk_1\\ \bdb\in\bk_{-1}\endSb
q^{\bda_1n_1+n_2(m_1+\bda_1)+ bdb_2m_2+\bdb_2(m_1+\bda_1)}
Q_{\bda+(3m_1,3n_1)+\bdb+(3m_2-1,3n_2-1)}P_{\bdb}P_{\bda}
\endalign$$
\noindent (the second term and the fourth term are negative to
each other)
$$\align
=&-q^{m_1n_2}(-\mu
q^{-(m_1+m_2)(n_1+n_2)}P_{-(3(m_1+m_2)-1,3(n_1+n_2)-1)}\\
&-\sum_{\bda,\bbda\in\bk_1}
q^{\bda(n_1+n_2)+\bbda_2(m_1+m_2)+\bda_1\bbda_2}Q_{(3(m_1+m_2)-1,3(n_1+n_2)-1)+\bda+\bbda}P_{\bbda}P_{\bda}\\
&-\sum\Sb \bda\in\bk_1\\
\bdb\in\bk_{-1}\endSb
q^{\bda_1(n_1+n_2)+\bdb_2(m_1+m_2)+\bdb_2+\bda_1}Q_{\bda+\bdb+(3(m_1+m_2)-1,3(n_1+n_2)-1)}P_{\bdb}P_{\bda})\\
=&-q^{m_1n_2}\eu_{21}(m_1+m_2,n_1+n_2).
\endalign$$

The following seven brackets are easily verified.
$$\align  & [\eu_{11}(m_1,n_1),\eu_{22}(m_2,n_2)]=0\\
& [\eu_{11}(m_1,n_1),\eu_{23}(m_2,n_2)]=0\\
 &[\eu_{11}(m_1,n_1),\eu_{31}(m_2,n_2)]=-q^{m_1n_2}\eu_{31}(m_1+m_2,n_1+n_2)\\
& [\eu_{11}(m_1,n_1),\eu_{32}(m_2,n_2)]=0\\
 & [\eu_{11}(m_1,n_1),\eu_{33}(m_2,n_2)]=0\\
 & [\eu_{12}(m_1,n_1),\eu_{12}(m_2,n_2)]=0\\
 & [\eu_{12}(m_1,n_1),\eu_{13}(m_2,n_2)]=0.
\endalign$$

Next, we have
$$\align  & [\eu_{12}(m_1,n_1),\eu_{21}(m_2,n_2)]\\
=&\mu
q^{-m_2n_2}\delta_{(m_1,n_1),(-m_2,n_2)}\\
&+\sum_{\bda'\in\bk_1}q^{n_2\bda'_1+n_1m_2+n_1\bda'_1}
Q_{(3m_1+1,3n_1+1)+\bda'+(3m_2-1,3n_2-1)P_{\bda'}}\\
&+\sum_{\bda\in\bk_1}q^{n_2m_1+\bda_2m_2+\bda_2m_1}
Q_{\bda+(3m_1+1,3n_1+1)+(3m_2-1,3n_2-1)}P_{\bda}\\
&+\sum_{\bdb\in\bk_{-1}}
q^{n_2m_1+\bdb_2m_2+\bdb_2m_1}Q_{(3m_1+1,3n_1+1)+\bdb+(3m_2-1,3n_2-1)}P_{\bdb}\endalign
$$
$$\align  =&q^{n_1m_2} (\sum_{\bda'\in\bk_1}q^{(n_1+n_2)\bda'_1}
Q_{(3m_1+3m_2,3n_1+3n_2)+\bda'}P_{\bda'}+\frac{1}{2}\delta_{(m_1,n_1),(-m_2,n_2)})\\
&-q^{n_2m_1}(-\sum_{\bda\in\bk_1}q^{\bda_2(m_2+m_1)}
Q_{\bda+(3m_1+3m_2,3n_1+3n_2)}P_{\bda}\\
&-\sum_{\bdb\in\bk_{-1}}
q^{\bdb_2(m_2+m_1)}Q_{(3m_1+3m_2,3n_1+3n_2)+\bdb}P_{\bdb}-\frac{1}{2}\delta_{(m_1,n_1),(-m_2,n_2)})\\
=&q^{n_1m_2}\eu_{11}(m_1+m_2,n_1+n_2)-q^{n_2m_1}\eu_{22}(m_1+m_2,n_1+n_2).
\endalign$$

The following six brackets can be checked easily.
$$\align  &
[\eu_{12}(m_1,n_1),\eu_{22}(m_2,n_2)]=q^{n_1m_2}\eu_{12}(m_1+m_2,n_1+n_2)\\
 &
 [\eu_{12}(m_1,n_1),\eu_{23}(m_2,n_2)]=q^{n_1m_2}\eu_{13}(m_1+m_2,n_1+n_2)\\
 &
 [\eu_{12}(m_1,n_1),\eu_{31}(m_2,n_2)]=-q^{m_1n_2}\eu_{32}(m_1+m_2,n_1+n_2)\\
& [\eu_{12}(m_1,n_1),\eu_{32}(m_2,n_2)]=0\\
 & [\eu_{12}(m_1,n_1),\eu_{33}(m_2,n_2)]=0\\
  & [\eu_{13}(m_1,n_1),\eu_{13}(m_2,n_2)]=0.
\endalign$$
$$\align  & [\eu_{13}(m_1,n_1),\eu_{21}(m_2,n_2)]\\
=&[\sum_{\bdb\in\bk_{-1}} q^{\bdb_1n_1}
Q_{(3m_1+2,3n_1+2)+\bdb}P_{\bdb}, -q^{-m_2n_2}\mu
P_{-(3m_2-1,3n_2-1)}\\
&-\sum_{\bda,bda'\in\bk_1}
q^{n_2\bda'_1+\bda_2m_2+\bda_2\bda'_1}Q_{\bda+\bda'+(3m_2-1,3n_2-1)}P_{\bda}P_{\bda'}\\
&-\sum\Sb \bda\in\bk_1\\ \bdb\in\bk_{-1}\endSb
q^{n_2\bda_1+\bdb_2m_2+\bdb_2\bda_1}
Q_{\bda+\bdb+(3m_2-1,3n_2-1)}P_{\bda}P_{\bdb}]\endalign$$
$$\align
=&-\sum\Sb \bda\in\bk_1\\ \bdb\in\bk_{-1}\endSb
q^{(\bda_1+\bdb_1+m_2)n_1+n_2+\bda_1+\bdb_2m_2+\bdb_2\bda_1}
Q_{(3m_1+2,3n_1+2)+\bda+\bdb+(3m_2-1,3n_2-1)}P_{\bda}P_{\bdb}\\
&+q^{-m_2n_2+(-m_1-m_2)n_1}\mu P_{(-3m_1-3m_2-1,-3n_1-3n_2-1)}\\
&+\sum\Sb \bda\in\bk_1\\ \bdb\in\bk_{-1}\endSb
q^{n_2(m_1+\bdb_1)+\bda_2m_2+\bda_2(m_1+\bdb_1)+\bdb_1n_1}
Q_{\bda+(3m_1+2,3n_1+2)+\bdb+(3m_2-1,3n_2-1)}P_{\bda}P_{\bdb}\\
&+\sum\Sb\bda'\in\bk_1\\ \bdb\in\bk_{-1}\endSb
q^{n_2\bda'_1+(n_1+\bdb_2)m_2+(n_1+\bdb_2)\bda'_1+\bdb_1n_1}
Q_{\bda'+(3m_1+2,3n_1+2)+\bdb+(3m_2-1,3n_2-1)}P_{\bda'}P_{\bdb}\\
&+\sum_{\bdb,\bdb'\in\bk_{-1}}
q^{n_2(\bdb'_1+m_1)+\bdb_2m_2+\bdb_2(\bdb'_1+m_1)+\bdb'_1n_1}
Q_{\bdb'+(3m_1+2,3n_1+2)+\bdb+(3m_2-1,3n_2-1)}P_{\bdb'}P_{\bdb}
\endalign$$
\noindent (the first and the fourth terms cancel each other)
$$\align =&q^{n_2m_1}(q^{-(m_1+m_2)(n_1+n_2)}\mu
P_{-(3m_1+3m_2-1,3n_1+3n_2-1)}\\
&+\sum_{\bda,bda'\in\bk_1}
q^{(n_1+n_2)\bda'_1+\bda_2(m_1+m_2)+\bda_2\bda'_1}Q_{\bda+\bda'+(3m_1+3m_2-1,3n_1+3n_2-1)}P_{\bda}P_{\bda'}\\
&-\sum\Sb \bda\in\bk_1\\ \bdb\in\bk_{-1}\endSb
q^{(n_1+n_2)\bda_1+\bdb_2(m_1+m_2)+\bdb_2\bda_1}
Q_{\bda+\bdb+(3m_1+3m_2-1,3n_1+3n_2-1)}P_{\bda}P_{\bdb})\\
=&-q^{n_2m_1}\eu_{23}(m_1+m_2,n_1+n_2).
\endalign$$

The following two brackets are easy.
$$\align  & [\eu_{13}(m_1,n_1),\eu_{22}(m_2,n_2)]=0\\
 & [\eu_{13}(m_1,n_1),\eu_{23}(m_2,n_2)]=0.
\endalign$$
$$\align  &-[\eu_{13}(m_2,n_2),\eu_{31}(m_1,n_1)]= [\eu_{31}(m_1,n_1),\eu_{13}(m_2,n_2)]\\
=&\sum\Sb \bda\in\bk_1\\\bdb\in\bk_{-1}\endSb
q^{n_1\bda_1+n_2\bdb_1} [Q_{(3m_1-2,3n_1-2)+\bda}P_{\bda},
Q_{(3m_2+2,3n_2+2)+\bdb}P_{\bdb}]\\
=&\sum_{\bdb\in\bk_{-1}} q^{n_2\bdb_1+n_1(m_2+\bdb_1)}
Q_{(3m_1+3m_2,3n_1+3n_2)+\bdb}P_{\bdb}\\
&-\sum_{\bda\in\bk_1} q^{n_1\bda_1+n_2(m_1+\bda_1)} Q_{(3m_1+3m_2,
3n_1+3n_2)+\bda} P_{\bda}\endalign$$
$$\align
=&q^{n_1m_2}(\sum_{\bdb\in\bk_{-1}} q^{n_2\bdb_1+n_1(m_2+\bdb_1)}
Q_{(3m_1+3m_2,3n_1+3n_2)+\bdb}P_{\bdb}+\frac{1}{2}\mu\delta_{(m_1+m_2,n_1+n_2),(0,0)})\\
&-q^{n_2m_1}(\sum_{\bda\in\bk_1} q^{n_1\bda_1+n_2(m_1+\bda_1)}
Q_{(3m_1+3m_2, 3n_1+3n_2)+\bda}
P_{\bda}+\frac{1}{2}\mu\delta_{(m_1+m_2,n_1+n_2),(0,0)})\\
=&q^{n_1m_2}\eu_{33}(m_1+m_2,n_1+n_2)-q^{n_2m_1}\eu_{11}(m_1+m_2,n_1+n_2).
\endalign$$

The following two brackets are easy.
$$\align  &
[\eu_{13}(m_1,n_1),\eu_{32}(m_2,n_2)]=q^{n_1m_2}\eu_{12}(m_1+m_2,n_1+n_2)\\
&[\eu_{13}(m_1,n_1),\eu_{33}(m_2,n_2)]=q^{n_1m_2}\eu_{13}(m_1+m_2,n_1+n_2).
\endalign$$
$$\align  & [\eu_{21}(m_1,n_1),\eu_{21}(m_2,n_2)]\\
=&\mu q^{-m_1n_1}\sum_{\bda\in\bk_1}
q^{(-m_1-m_2-\bda_1)n_2+\bda_2m_2+\bda_2(-m_1-m_2-\bda_1)}P_{\bda}P_{(-3m_1-3m_2,-3n_1-3n_2)-\bda}\\
&-\mu q^{-m_2n_2}\sum_{\bda'\in\bk_1}
q^{n_1\bda'_1+(-n_1-n_2-\bda'_2)m_1+(-n_1-n_2-\bda'_2)\bda'_1}P_{\bda'}P_{(-3m_1-3m_2,-3n_1-3n_2)-\bda'}\\
&+[\sum_{\bda,\bda'\in\bk_1}q^{n_1\bda'_1+\bda_2m_1+\bda_2\bda'_1}
Q_{\bda+\bda'+(3m_1-1,3n_1-1)}P_{\bda}P_{\bda'},\\
&\sum_{\bda,\bda'\in\bk_1}q^{n_2\bda'_1+\bda_2m_2+\bda_2\bda'_1}
Q_{\bda+\bda'+(3m_2-1,3n_2-1)}P_{\bda}P_{\bda'}]\\
&+[\sum_{\bda,\bda'\in\bk_1}q^{n_1\bda'_1+\bda_2m_1+\bda_2\bda'_1}
Q_{\bda+\bda'+(3m_1-1,3n_1-1)}P_{\bda}P_{\bda'},\\
&\sum\Sb \bda\in\bk_1\\\bdb\in\bk_{-1}\endSb
q^{n_2\bda_1+\bdb_2m_2+\bdb_2\bda_1}Q_{\bda+\bdb+(3m_2-1,3n_2-1)}P_{\bda}P_{\bdb}]\\
&+[\sum\Sb \bda\in\bk_1\\\bdb\in\bk_{-1}\endSb
q^{n_1\bda_1+\bdb_2m_1+\bdb_2\bda_1}Q_{\bda+\bdb+(3m_1-1,3n_1-1)}P_{\bda}P_{\bdb},\\
&\sum_{\bda,\bda'\in\bk_1}q^{n_2\bda'_1+\bda_2m_2+\bda_2\bda'_1}
Q_{\bda+\bda'+(3m_2-1,3n_2-1)}P_{\bda}P_{\bda'}]\\
&+[\sum\Sb \bda\in\bk_1\\\bdb\in\bk_{-1}\endSb
q^{n_1\bda_1+\bdb_2m_1+\bdb_2\bda_1}Q_{\bda+\bdb+(3m_1-1,3n_1-1)}P_{\bda}P_{\bdb},\\
&\sum\Sb \bda\in\bk_1\\\bdb\in\bk_{-1}\endSb
q^{n_2\bda_1+\bdb_2m_2+\bdb_2\bda_1}Q_{\bda+\bdb+(3m_2-1,3n_2-1)}P_{\bda}P_{\bdb}]
\endalign$$
\noindent (the first  and the second terms cancel each other)
$$\align =&\sum_{\bda,\bbda,\bbda'\in\bk_1}
q^{n_1(\bbda_1+\bbda'_1+m_2)+\bda_2m_1+\bda_2(\bbda_1+\bbda'_1+m_2)+n_2\bbda'_1+\bbda_2m_2+\bbda_2\bbda'_1}\\
&.Q_{\bda+\bbda+\bbda'+(3m_1-1,3n_1-1)+(3m_2-1,3n_2-1)}P_{\bda}P_{\bbda}P_{\bbda'}\\
&+\sum_{\bda',\bbda,\bbda'\in\bk_1}
q^{n_1\bda'_1+(\bbda_2+\bbda'_2+n_2)m_1+(\bbda_2+\bbda'_2+n_2)\bda'_1+n_2\bbda'_1+\bbda_2m_2+\bbda_2\bbda'_1}\\
&.Q_{\bda'+\bbda+\bbda'+(3m_1-1,3n_1-1)+(3m_2-1,3n_2-1)}P_{\bda'}P_{\bbda}P_{\bbda'}\\
\endalign$$
$$\align&-\sum_{\bda,\bda',\bbda\in\bk_1}
q^{n_2(\bda_1+\bda'_1+m_1)+\bbda_2m_2+\bbda_2(\bda_1+\bda'_1+m_1)+n_1\bda'_1+\bda_2m_1+\bda_2\bda'_1}\\
&.Q_{\bda+\bda'+\bbda+(3m_1-1,3n_1-1)+(3m_2-1,3n_2-1)}P_{\bda}P_{\bda'}P_{\bbda}\\
&-\sum_{\bda,\bda',\bbda'\in\bk_1}
q^{n_2\bbda'_1+(\bda_2+\bda'_2+n_1)m_2+(\bda_2+\bda'_2+n_1)\bbda'_1+n_1\bda'_1+\bda_2m_1+\bda_2\bda'_1}\\
&.Q_{\bda+\bda'+\bbda'+(3m_1-1,3n_1-1)+(3m_2-1,3n_2-1)}P_{\bda}P_{\bda'}P_{\bbda'}\\
 &-\sum\Sb \bda,\bda'\in\bk_1\\ \bdb\in\bk_{-1}\endSb
q^{n_2(\bda_1+\bda'_1+m_1)+\bdb_2m_2+\bdb_2(\bda_1+\bda'_1+m_1)+n_1\bda'_1+\bda_2m_1\bda_2\bda'_1}\\
&.Q_{\bda+\bda'+\bdb+(3m_1-1,3n_1-1)+(3m_2-1,3n_2-1)}P_{\bda}P_{\bda'}P_{\bdb}\\
&+\sum\Sb \bda,\bda'\in\bk_1\\ \bdb\in\bk_{-1}\endSb
q^{n_1(\bda_1+\bda'_1+m_2)+\bdb_2m_1+\bdb_2(\bda_1+\bda'_1+m_2)+n_2\bda'_1+\bda_2m_2+\bda_2\bda'_1}\\
&.Q_{\bda+\bda'+\bdb+(3m_1-1,3n_1-1)+(3m_2-1,3n_2-1)}P_{\bda}P_{\bda'}P_{\bdb}\\
&+\sum\Sb \bda,\bbda\in\bk_1\\ \bdb\in\bk_{-1}\endSb
q^{n_1\bda_1+(\bbda_2+\bdb_2+n_2)m_1+(\bbda_2+\bdb_2+n_2)\bda_1+n_2\bbda_1+\bdb_2m_2+\bdb_2\bbda_1}\\
&.Q_{\bda+\bbda+\bdb+(3m_1-1,3n_1-1)+(3m_2-1,3n_2-1)}P_{\bda}P_{\bbda}P_{\bdb}\\
&-\sum\Sb \bda,\bbda\in\bk_1\\ \bbdb\in\bk_{-1}\endSb
q^{n_2\bda_1+(\bbda_2+\bbdb_2+n_1)m_2+(\bbda_2+\bbdb_2+n_1)\bda_1+n_1\bbda_1+\bbdb_2m_1+\bbdb_2\bbda_1}\\
&.Q_{\bda+\bbda+\bbdb+(3m_1-1,3n_1-1)+(3m_2-1,3n_2-1)}P_{\bda}P_{\bbda}P_{\bbdb}
=0 \endalign$$ as the first term cancels the fourth, the second
term cancels the third, the fifth term cancels the seventh, and
the sixth term cancels the eighth.

$$\align  & -[\eu_{21}(m_2,n_2),\eu_{22}(m_1,n_1)] =[\eu_{22}(m_1,n_1),\eu_{21}(m_2,n_2)]\\
=&\mu\sum_{\bda\in\bk_1}q^{\bda_2m_1-m_2n_2}[Q_{(3m_1,3n_1)+\bda}P_{\bda},P_{-(3m_2-1,3n_2-1)}]\\
&+\sum_{\bda,\bda',\bbda\in\bk_1}q^{\bda_2m_1+n_2\bda'_1+\bbda_2m_2+\bbda_2\bda'_1}[Q_{(3m_1,3n_1)+\bda}P_{\bda},
Q_{\bbda+\bda'+(3m_2-1,3n_2-1)}P_{\bda'}P_{\bbda}]\\
&+\sum\Sb\bda,\bbda\in\bk_1\\ \bdb\in\bk_{-1}\endSb
q^{\bda_2m_1+n_2\bbda_1+\bdb_2m_2+\bdb_2\bbda_1}[Q_{(3m_1,3n_1)+\bda}P_{\bda},
Q_{\bbda+\bdb+(3m_2-1,3n_2-1)}P_{\bda}P_{\bdb}]\\
&+\sum\Sb \bda\in\bk_1\\
\bdb,\bbdb\in\bk_{-1}\endSb
q^{\bdb_2m_1+n_2\bda_1+\bbdb_2m_2+\bbdb_2\bda_1}[Q_{(3m_1,3n_1)+\bdb}P_{\bdb},
Q_{\bda+\bbdb+(3m_2-1,3n_2-1)}P_{\bda}P_{\bbdb}]\\
=&-\mu q^{-m_2n_2+(-n_1-n_2)m_1}P_{-3m_1-3m_2+1,-3n_1-3n_2+1}\\
&+\sum_{\bda,\bda'\in\bk_1}q^{(\bda_2+\bda'_2+n_2)m_1+n_2\bda'_1+\bda_2m_2+\bda_2\bda'_1}
Q_{(3m_1,3n_1)+\bda+\bda'+(3m_2-1,3n_2-1)}P_{\bda}P_{\bda'}\\
&-\sum_{\bda,\bbda\in\bk_1}q^{n_2(m_1+\bbda_1)+\bda_2m_2+\bda_2(m_1+\bbda_1)+\bbda_2m_1}
Q_{\bda+\bbda+(3m_1,3n_1)+(3m_2-1,3n_2-1)}P_{\bda}P_{\bbda}\\
&-\sum_{\bda',\bbda\in\bk_1}q^{n_2\bda'_1+(\bbda_2+n_1)m_2+(\bbda_2+n_1)\bda'_1+\bbda_2m_1}
Q_{\bda'+\bbda+(3m_1,3n_1)+(3m_2-1,3n_2-1)}P_{\bda'}P_{\bbda}\\
&-\sum\Sb \bda\in\bk_1\\ \bdb\in\bk_{-1}\endSb
q^{n_2(\bda_1+m_1)+\bdb_2m_2+\bdb_2(\bda_1+m_1)+\bda_2m_1}Q_{\bda+\bdb+(3m_1,3n_1)+(3m_2-1,3n_2-1)}P_{\bda}P_{\bdb}\\
&+\sum\Sb \bda\in\bk_1\\ \bdb\in\bk_{-1}\endSb
q^{(\bda_2+\bdb_2+n_2)m_1+n_2\bda_1+\bdb_2m_2+\bdb_2\bda_1}Q_{\bda+\bdb+(3m_1,3n_1)+(3m_2-1,3n_2-1)}P_{\bda}P_{\bdb}\\
&-\sum\Sb \bda\in\bk_1\\ \bdb\in\bk_{-1}\endSb
q^{n_2\bda_1+(n_1+\bdb_2)m_2+(n_1+\bdb_2)\bda_1+\bdb_2m_1}Q_{\bda+\bdb+(3m_1,3n_1)+(3m_2-1,3n_2-1)}P_{\bda}P_{\bdb}
\endalign$$
\noindent(the second and the third terms are cancelled,  and the
fifth  and the sixth terms are cancelled)
$$\align
=&q^{n_1m_2}(-\mu
q^{-(n_1+n_2)(m_1+m_2)}P_{-(3(m_1+m_2)-1,3(n_1+n_2)-1)}\\
&-\sum_{\bbda,\bda'\in\bk_1}q^{(n_1+n_2)\bda'_1+(m_1+m_2)\bbda_2+\bbda_2\bda'_1}
Q_{(3m_1+3m_2-1,3n_1+3n_2-1)+\bda'+\bbda}P_{\bda'}P_{\bbda}\\
&-\sum\Sb \bda\in\bk_1\\
\bdb\in\bk_{-1}\endSb
q^{(n_1+n_2)\bda_1+(m_1+m_2)\bdb_2+\bdb_2\bda_1})
Q_{(3m_1+3m_2-1,3n_1+3n_2-1)+\bda+\bdb}P_{\bda}P_{\bdb}\\
=&q^{n_1m_2}\eu_{21}(m_1+m_2,n_1+n_2).
\endalign$$
$$\align  & [\eu_{21}(m_1,n_1),\eu_{23}(m_2,n_2)]\\
=&\mu q^{-m_1n_1}\sum_{\bda\in\bk_1}
q^{n_2(-m_1-m_2-\bda_1)+\bda_2m_2+\bda_2(-m_1-m_2-\bda_1)}
P_{\bda}P_{(-3m_1-3m_2,-3n_1-3n_2)-\bda}\\
&+[\sum_{\bda,\bda'\in\bk_1}q^{n_1\bda'_1+\bda_2m_1+\bda_2\bda'_1}Q_{\bda+\bda'+(3m_1-1,3n_1-1)}P_{\bda}P_{\bda'},\\
&\sum\Sb \bda\in\bk_1\\ \bdb\in\bk_{-1}\endSb
q^{n_2\bdb_1+\bda_2m_2+\bda_2\bdb_1}
Q_{(3m_2+1,3n_2+1)+\bda+\bdb}P_{\bda}P_{\bdb}]\\
&-\mu q^{-m_2n_2} \sum_{\bda\in\bk_1}
q^{n_1\bda_1+(-n_1-n_2-\bda_2)m_1+(-n_1-n_2-\bda_2)\bda_1}
P_{\bda}P_{(-3m_1-3m_2,-3n_1-3n_2)-\bda}\\
&+[\sum\Sb \bda\in\bk_1\\ \bdb\in\bk_{-1}\endSb
q^{n_1\bda_1+\bdb_2m_1+\bdb_2\bda_1}
Q_{\bda+\bdb+(3m_1-1,3n_1-1)}P_{\bda}P_{\bdb},\\
&\sum\Sb \bda\in\bk_1\\ \bdb\in\bk_{-1}\endSb
q^{n_2\bdb_1+\bda_2m_2+\bda_2\bdb_1}Q_{(3m_2+1,3n_2+1)}P_{\bda}P_{\bdb}]\\
&+[\sum\Sb \bda\in\bk_1\\ \bdb\in\bk_{-1}\endSb
q^{n_1\bda_1+\bdb_2m_1+\bdb_2\bda_1}
Q_{\bda+\bdb+(3m_1-1,3n_1-1)}P_{\bda}P_{\bdb},\\
&\sum_{\bdb,\bdb'\in\bk_{-1}}q^{n_2\bdb'_1+\bdb_2m_2+\bdb_2\bdb'_1}
Q_{(3m_2+1,3n_2+1)+\bdb+\bdb'}P_{\bdb}P_{\bdb'}]
\endalign$$
\noindent (the first term and the third term are negative to each
other)
$$\align  =&\sum\Sb
\bda,\bbda\in\bk_1\\\bdb\in\bk_{-1}\endSb
q^{n_1(m_2+\bbda_1+\bdb_1)+\bda_2m_1+\bda_2(m_1+\bbda_1+\bdb_1)+n_2\bdb_1+\bbda_2m_2+\bbda_2\bdb_1}\\
&.Q_{\bda+\bbda+\bdb+(3m_1+3m_2,3n_1+3n_2)}P_{\bda}P_{\bbda}P_{\bdb}\\
&+\sum\Sb \bda',\bbda\in\bk_1\\\bdb\in\bk_{-1}\endSb
q^{n_1\bda'_1+(n_2+\bbda_2+\bdb_2)m_1+(n_2+\bbda_2+\bdb_2)\bbda_1+n_2\bdb_1+\bbda_2m_2+\bbda_2\bdb_1}\\
&.Q_{\bda'+\bbda+\bdb+(3m_1+3m_2,3n_1+3n_2)}P_{\bda'}P_{\bbda}P_{\bdb}\\
&-\sum\Sb \bda,\bda'\in\bk_1\\\bdb\in\bk_{-1}\endSb
q^{n_2\bdb_1+(\bda_2+\bda'_2+n_1)m_2+(\bda_2+\bda'_2+n_1)\bdb_1+n_1\bda'_1+\bda_2m_1+\bda_2\bda'_1}\\
&.Q_{\bda'+\bda+\bdb+(3m_1+3m_2,3n_1+3n_2)}P_{\bda'}P_{\bda}P_{\bdb}\\
\endalign$$
$$\align&+\sum\Sb \bda\in\bk_1\\\bdb,\bdb'\in\bk_{-1}\endSb
q^{n_1(m_2+\bda_1+\bdb'_1)+\bdb_2m_1+\bdb_2(m_2+\bda_1+\bdb'_1)+n_2\bdb'_1+\bda_2m_2+\bda_2\bdb'_1}\\
&.Q_{\bda+\bdb+\bdb'+(3m_1+3m_2,3n_1+3n_2)}P_{\bda}P_{\bdb}P_{\bdb'}\\
 &-\sum\Sb \bda,\bbda\in\bk_1\\\bdb\in\bk_{-1}\endSb
q^{n_2(\bda_1+\bdb_1+m_1)+\bbda_2m_2+\bbda_2(\bda_1+\bdb+m_1)+n_1\bda_1+\bdb_2m_1+\bdb_2\bda_1}\\
&.Q_{\bda+\bbda+\bdb+(3m_1+3m_2,3n_1+3n_2)}P_{\bda}P_{\bbda}P_{\bdb}\\
&+\sum\Sb \bda\in\bk_1\\\bdb,\bdb'\in\bk_{-1}\endSb
q^{n_1\bda_1+(n_2+\bdb_2+\bdb'_2)m_1+(n_2+\bdb_2+\bdb'_2)\bda_1+n_2\bdb'_1+\bdb_2m_2+\bdb_2\bdb'_1}\\
&.Q_{\bda+\bdb+\bdb'+(3m_1+3m_2,3n_1+3n_2)}P_{\bda}P_{\bdb}P_{\bdb'}\\
&-\sum\Sb \bda\in\bk_1\\\bdb,\bbdb\in\bk_{-1}\endSb
q^{n_2(\bda_1+\bbdb_1+m_1)+\bdb_2m_2+\bdb_2(\bda_1+\bbdb_1+m_1)+n_1\bda_1+\bbdb_2m_1+\bbdb_2\bda_1}\\
&.Q_{\bda+\bdb+\bbdb+(3m_1+3m_2,3n_1+3n_2)}P_{\bda}P_{\bdb}P_{\bbdb}\\
&-\sum\Sb \bda\in\bk_1\\\bdb',\bbdb\in\bk_{-1}\endSb
q^{n_2\bdb'_1+(\bda_2+\bbdb_2+n_1)m_2+(\bda_2+\bbdb_2+n_1)\bdb'_1+n_1\bda_1+\bbdb_2m_1+\bbdb_2\bda_1}\\
&.Q_{\bda+\bdb'+\bbdb+(3m_1+3m_2,3n_1+3n_2)}P_{\bda}P_{\bdb'}P_{\bbdb}=0
\endalign$$
as the first term and the third, the second and the fifth, the
fourth and the eighth, the sixth and the seventh are cancelled.

The following three brackets are easy.
$$\align  & [\eu_{21}(m_1,n_1),\eu_{31}(m_2,n_2)]=0 \\
 &
 [\eu_{21}(m_1,n_1),\eu_{32}(m_2,n_2)]=-q^{n_2m_1}\eu_{31}(m_1+m_2,n_1+n_2)\\
  & [\eu_{21}(m_1,n_1),\eu_{33}(m_2,n_2)]=0.
\endalign$$
$$\align  & [\eu_{22}(m_1,n_1),\eu_{22}(m_2,n_2)]\\
=&\sum_{\bda,\bda'\in\bk_1}
q^{\bda_2m_1+\bda'_2m_2}[Q_{(3m_1,3n_1)+\bda}P_{\bda},Q_{(3m_2,3n_2)+\bda'}P_{\bda'}]\\
&+\sum_{\bdb,\bdb'\in\bk_{-1}}
q^{\bdb_2m_1+\bdb'_2m_2}[Q_{(3m_1,3n_1)+\bdb}P_{\bdb},Q_{(3m_2,3n_2)+\bdb'}P_{\bdb'}]\\
=&\sum_{\bda'\in\bk_1} q^{(n_2+\bda'_2)m_1+\bda'_2m_2}
Q_{(3m_1,3n_1)+(3m_2,3n_2)+\bda'}P_{\bda'}\\
&-\sum_{\bda\in\bk_1} q^{\bda_2m_1+(n_1+\bda_2)m_2}
Q_{(3m_1,3n_1)+(3m_2,3n_2)+\bda}P_{\bda}\\
&+\sum_{\bdb'\in\bk_{-1}} q^{(n_2+\bdb'_2)m_1+\bdb'_2m_2}
Q_{(3m_1,3n_1)+(3m_2,3n_2)+\bdb'}P_{\bdb'}\\
&-\sum_{\bdb\in\bk_{-1}} q^{\bdb_2m_1+(n_1+\bdb_2)m_2}
Q_{(3m_1,3n_1)+(3m_2,3n_2)+\bdb}P_{\bdb}\\
=&q^{n_1m_2}(-\sum_{\bda\in\bk_1}
q^{\bda_2(m_1+m_2)}Q_{(3m_1,3n_1)+(3m_2,3n_2)+\bda}P_{\bda}\\
&-\sum_{\bdb\in\bk_{-1}}
q^{\bdb_2(m_1+m_2)}Q_{(3m_1,3n_1)+(3m_2,3n_2)+\bdb}P_{\bdb}-\frac{1}{2}\mu\delta_{(m_1+m_2,n_1+n_2),(0,0)})\\
&-q^{n_2m_1}(-\sum_{\bda\in\bk_1}
q^{\bda_2(m_1+m_2)}Q_{(3m_1,3n_1)+(3m_2,3n_2)+\bda}P_{\bda}\\
&-\sum_{\bdb\in\bk_{-1}}
q^{\bdb_2(m_1+m_2)}Q_{(3m_1,3n_1)+(3m_2,3n_2)+\bdb}P_{\bdb}-\frac{1}{2}\mu\delta_{(m_1+m_2,n_1+n_2),(0,0)})\\
=&q^{n_1m_2}\eu_{22}(m_1+m_2,n_1+n_2)-q^{n_2m_1}\eu_{22}(m_1+m_2,n_1+n_2).
\endalign$$

$$\align  & [\eu_{22}(m_1,n_1),\eu_{23}(m_2,n_2)]\\
=&[-\sum_{\bda\in\bk_1} q^{\bda_2m_1}
Q_{(3m_1,3n_1)+\bda}P_{\bda}-\sum_{\bdb\in\bk_{-1}}q^{\bdb_2m_1}Q_{(3m_1,3n_1)+\bdb}P_{\bdb},\\
&-q^{-m_2n_2}\mu P_{-(3m_2+1,3n_2+1)}-\sum\Sb \bda\in\bk_1\\
\bdb\in\bk_{-1}\endSb q^{n_2\bdb_1+\bda_2m_2+\bda_2\bdb_1}
Q_{\bda+\bdb+(3m_2+1,3n_2+1)}P_{\bda}P_{\bdb}\\
&-\sum_{\bdb,\bdb'\in\bk_{-1}}
q^{n_2\bdb'_1+\bdb_2m_2+\bdb_2\bdb'_1}
Q_{\bdb'+\bdb+(3m_2+1,3n_2+1)}P_{\bdb}P_{\bdb'}]\\
=&\sum\Sb \bda\in\bk_1\\ \bdb\in\bk_{-1}\endSb
q^{(n_2+\bda_2+\bdb_2)m_1+n_2\bdb_1+\bda_2m_2+\bda_2\bdb_1}
Q_{(3m_1,3n_1)+\bda+\bdb+(3m_2+1,3n_2+1)} P_{\bda}P_{\bdb}\\
&-\sum\Sb \bda\in\bk_1\\ \bdb\in\bk_{-1}\endSb
q^{n_2\bdb_1+(n_1+\bda_2)m_2+(\bda_2+n_1)\bdb_1+\bda_2m_1}
 Q_{(3m_1,3n_1)+\bda+\bdb+(3m_2+1,3n_2+1)} P_{\bda}P_{\bdb}\\
 &-q^{-m_2n_2+(-n_2-n_1)m_1} \mu
 P_{-(3m_2+1,3n_2+1)-(3m_1,3n_1)}\\
 &-\sum\Sb \bda\in\bk_1\\ \bdb\in\bk_{-1}\endSb
 q^{n_2(m_1+\bdb_1)+\bda_2m_2+\bda_2(m_1+\bdb_1)+\bdb_2m_1}
  Q_{(3m_1,3n_1)+\bda+\bdb+(3m_2+1,3n_2+1)} P_{\bda}P_{\bdb}\\
  &+\sum_{\bdb,\bdb'\in\bk_{-1}}
  q^{(n_2+\bdb_2+\bdb'_2)m_1+n_2\bdb'_1+\bdb_2m_2+\bdb_2\bdb'_1}
  Q_{(3m_1,3n_1)+\bdb+\bdb'+(3m_2+1,3n_2+1)} P_{\bdb}P_{\bdb'}\\
&-\sum_{\bdb,\bdb'\in\bk_{-1}}
q^{n_2(\bdb'_1+m_1)+\bdb_2m_2+\bdb_2(\bdb'_1+m_1)+\bdb'_2m_1}
Q_{(3m_1,3n_1)+\bdb+\bdb'+(3m_2+1,3n_2+1)} P_{\bdb}P_{\bdb'}\\
&-\sum_{\bdb,\bdb'\in\bk_{-1}}
q^{n_2\bdb'_1+(n_1+\bdb_2)m_2+(\bdb_2+n_1)\bdb'_1+\bdb_2m_1}
Q_{(3m_1,3n_1)+\bdb+\bdb'+(3m_2+1,3n_2+1)} P_{\bdb}P_{\bdb'}
\endalign$$
\noindent (the first term cancels the fourth while the fifth
cancels the sixth)
$$\align
=&q^{m_2n_1}(-q^{-(m_1+m_2)(n_1+n_2)}\mu
P_{-(3m_1+3m_2+1,3n_1+3n_2+1)}\\
&-\sum\Sb \bda\in\bk_1\\ \bdb\in\bk_{-1}\endSb
q^{\bda_2(m_1+m_2)+(n_1+n_2)\bdb_1+\bda_2\bdb_1}
Q_{\bda+\bdb+(3m_1+3m_2+1,3n_1+3n_2+1)}P_{\bda}P_{\bdb}\\
&-\sum_{\bdb,\bdb'\in\bk_{-1}}
q^{\bdb'_1(n_1+n_2)+\bdb_2(m_1+m_2)\bdb_2\bdb'_1})
Q_{\bdb+\bdb'+(3m_1+3m_2+1,3n_1+3n_2+1)}P_{\bdb}P_{\bdb'}\\
=&q^{m_2n_1}\eu_{23}(m_1+m_2,n_1+n_2).
\endalign$$

The following three brackets are easy.
$$\align  & [\eu_{22}(m_1,n_1),\eu_{31}(m_2,n_2)]=0\\
&
[\eu_{22}(m_1,n_1),\eu_{32}(m_2,n_2)]=-q^{n_2m_1}\eu_{32}(m_1+m_2,n_1+n_2)\\
 & [\eu_{22}(m_1,n_1),\eu_{33}(m_2,n_2)]=0.
\endalign$$

$$\align  & [\eu_{23}(m_1,n_1),\eu_{23}(m_2,n_2)]\\
=&[-q^{-m_1n_1}\mu P_{-(3m_1+1,3n_1+1)}-\sum\Sb \bda\in\bk_1\\
\bdb\in\bk_{-1}\endSb q^{n_1\bdb_1+\bda_2m_1+\bda_2\bdb_1}
Q_{\bda+\bdb+(3m_1+1,3n_1+1)}P_{\bda}P_{\bdb}\\
&-\sum_{\bdb,\bdb'\in\bk_{-1}}
q^{n_1\bdb'_1+\bdb_2m_1+\bdb_2\bdb'_1}
Q_{\bdb+\bdb'+(3m_1+1,3n_1+1)}P_{\bdb}P_{\bdb'},\\
&-q^{-m_2n_2}\mu P_{-(3m_2+1,3n_2+1)}-\sum\Sb \bda\in\bk_1\\
\bdb\in\bk_{-1}\endSb q^{n_2\bdb_1+\bda_2m_2+\bda_2\bdb_1}
Q_{\bda+\bdb+(3m_2+1,3n_2+1)}P_{\bda}P_{\bdb}\\
&-\sum_{\bdb,\bdb'\in\bk_{-1}}
q^{n_2\bdb'_1+\bdb_2m_2+\bdb_2\bdb'_1}
Q_{\bdb+\bdb'+(3m_2+1,3n_2+1)}P_{\bdb}P_{\bdb'}]\\
=&\sum_{\bdb\in\bk_{-1}}\mu
q^{-m_1n_1+n_2(-m_1-m_2-\bdb_1)+\bdb_2m_2+\bdb_2(-m_1-m_2-\bdb_1)}P_{\bdb}P_{-(3m_1+1,3n_1+1)-(3m_2+1,3n_2+1)-\bdb}\\
&-\sum_{\bdb\in\bk_{-1}}\mu
q^{-m_2n_2+n_1\bdb_1+(-n_2-n_1-\bdb_2)m_1+(-n_2-n_1-\bdb_2)\bdb_1}P_{\bdb}P_{-(3m_1+1,3n_1+1)-(3m_2+1,3n_2+1)-\bdb}\\
&+\sum_{\bdb,\bbdb,\bdb'\in\bk_{-1}}
q^{n_1(\bdb'_1+\bbdb_1+m_2)+\bdb_2m_1+\bdb_2(\bdb'_1+\bbdb_1+m_2)+n_2\bdb'_1+\bbdb_2m_2+\bbdb_2\bdb'_1}\\
&.Q_{\bdb+\bdb'+\bbdb+(3m_1+1,3n_1+1)+(3m_2+1,3n_2+1)} P_{\bdb}
P_{\bdb'} P_{\bbdb}\\
&+\sum_{\bdb,\bbdb,\bdb'\in\bk_{-1}}
q^{n_1\bdb'_1+(\bdb_2+\bbdb_2+n_2)m_1+\bdb'_1(\bdb_2+\bbdb_2+n_2)+n_2\bbdb_1+\bdb_2m_2+\bdb_2\bbdb_1}\\
&.Q_{\bdb+\bdb'+\bbdb+(3m_1+1,3n_1+1)+(3m_2+1,3n_2+1)} P_{\bdb}
P_{\bdb'} P_{\bbdb}\\
 &-\sum_{\bdb,\bbdb,\bdb'\in\bk_{-1}}
q^{n_2(\bdb'_1+\bbdb_1+m_1)+\bdb_2m_2+\bdb_2(\bdb'_1+\bbdb_1+m_1)+n_1\bdb'_1+\bbdb_2m_1+\bbdb_2\bdb'_1}\\
&.Q_{\bdb+\bdb'+\bbdb+(3m_1+1,3n_1+1)+(3m_2+1,3n_2+1)} P_{\bdb}
P_{\bdb'} P_{\bbdb}\endalign$$
$$\align
&-\sum_{\bdb,\bbdb,\bdb'\in\bk_{-1}}
q^{n_2\bdb'_1+(\bdb_2+\bbdb_2+n_1)m_2+(\bdb_2+\bbdb_2+n_1)\bdb'_1+n_1\bbdb_1+\bdb_2m_1+\bdb_2\bbdb_1}\\
&.Q_{\bdb+\bdb'+\bbdb+(3m_1+1,3n_1+1)+(3m_2+1,3n_2+1)} P_{\bdb}
P_{\bdb'} P_{\bbdb}\\
&+\sum\Sb \bda\in\bk_1\\ \bdb,\bdb'\in\bk_{-1}\endSb
q^{n_1(\bdb_1+\bdb'_1+m_2)+\bda_2m_1+\bda_2(\bdb_1+\bdb'_1+m_2)+n_2\bdb'_1+\bdb_2m_2+\bdb_2\bdb'_1}\\
&.Q_{\bda+\bdb+\bdb'+(3m_1+1,3n_1+1)+(3m_2+1,3n_2+1)} P_{\bda}
P_{\bdb} P_{\bdb'}\\
&+\sum\Sb \bda\in\bk_1\\ \bdb,\bdb'\in\bk_{-1}\endSb
q^{n_1\bdb_1+(\bda_2+\bdb'_2+n_2)m_1+(\bda_2+\bdb'_2+n_2)\bdb_1+n_2\bdb'_1+\bda_2m_2+\bda_2\bdb'_1}\\
&.Q_{\bda+\bdb+\bdb'+(3m_1+1,3n_1+1)+(3m_2+1,3n_2+1)} P_{\bda}
P_{\bdb} P_{\bdb'}\\
 &-\sum\Sb \bda\in\bk_1\\
\bdb,\bdb'\in\bk_{-1}\endSb
q^{n_2(\bdb_1+\bdb'_1+m_1)+\bda_2m_2+\bda_2(\bdb_1+\bdb'_1+m_1)+n_1\bdb'_1+\bdb_2m_1+\bdb_2\bdb'_1}\\
&.Q_{\bda+\bdb+\bdb'+(3m_1+1,3n_1+1)+(3m_2+1,3n_2+1)} P_{\bda}
P_{\bdb} P_{\bdb'}\\
&-\sum\Sb \bda\in\bk_1\\ \bdb,\bdb'\in\bk_{-1}\endSb
q^{n_2\bdb_1+(\bda_2+\bdb'_2+n_1)m_2+(\bda_2+\bdb'_2+n_1)\bdb_1+n_1\bdb'_1+\bda_2m_1+\bda_2\bdb'_1}\\
&.Q_{\bda+\bdb+\bdb'+(3m_1+1,3n_1+1)+(3m_2+1,3n_2+1)} P_{\bda}
P_{\bdb} P_{\bdb'}=0 \endalign$$  as the first two terms, the
third and the sixth, the fourth  and the fifth, the seventh and
the last term, are negative to each other.

$$\align  & [\eu_{23}(m_1,n_1),\eu_{31}(m_2,n_2)]\\
=&[-q^{-m_1n_1}\mu P_{-(3m_1+1,3n_1+1)}-\sum\Sb \bda\in\bk_1\\
\bdb\in\bk_{-1}\endSb q^{n_1\bdb_1+\bda_2m_1+\bda_2\bdb_1}
Q_{\bda+\bdb+(3m_1+1,3n_1+1)}P_{\bda}P_{\bdb}\\
&-\sum_{\bdb,\bdb'\in\bk_{-1}}
q^{n_1\bdb'_1+\bdb_2m_1+\bdb_2\bdb'_1}
Q_{\bdb+\bdb'+(3m_1+1,3n_1+1)}P_{\bdb}P_{\bdb'},\\
&\sum_{\bda\in\bk_1} q^{\bda_1n_2}
Q_{(3m_2-2,3n_2-2)+\bda}P_{\bda}]\endalign$$
$$\align
=&-q^{-m_1n_1+(-m_1-m_2)n_2}\mu
P_{-(3m_1+1,3n_1+1)-(3m_2-2,3n_2-2)}\\
&-\sum_{\bda,\bda'\in\bk_1}
q^{n_1(m_2+\bda'_1)+\bda_2m_1+\bda_2(m_2+\bda'_1)+\bda'_1n_2}
Q_{\bda+\bda'+(3m_2-2,3n_2-2)+(3m_1+1,3n_1+1)}P_{\bda}P_{\bda'}\\
&+\sum\Sb \bda\in\bk_1\\ \bdb\in\bk_{-1}\endSb
q^{(\bda_1+\bdb_1+m_1)n_2+\bda_2m_1+\bda_2\bdb_1+n_1\bdb_1}
Q_{\bda+\bdb+(3m_2-2,3n_2-2)+(3m_1+1,3n_1+1)}P_{\bda}P_{\bdb}\\
&-\sum\Sb \bda\in\bk_1\\ \bdb\in\bk_{-1}\endSb
q^{n_1(m_2+\bda_1)+\bdb_2m_1+\bdb_2(m_2+\bda_1)+\bda_1n_2}
Q_{\bda+\bdb+(3m_2-2,3n_2-2)+(3m_1+1,3n_1+1)}P_{\bda}P_{\bdb}\\
&-\sum\Sb \bda\in\bk_1\\ \bdb\in\bk_{-1}\endSb
q^{n_1\bdb_1+(\bda_2+n_2)m_1+(\bda_2+n_2)\bdb_1+\bda_1n_2}
Q_{\bda+\bdb+(3m_2-2,3n_2-2)+(3m_1+1,3n_1+1)}P_{\bda}P_{\bdb}
\endalign$$
\noindent (the third term and the fifth are cancelled)
$$\align
=&q^{n_1m_2}(-q^{-(m_1+m_2)(n_1+n_2)} \mu
P_{-(3m_1+3m_2-1,3n_1+3n_2-1)}\\
&-\sum_{\bda,\bda'\in\bk_1}
q^{(n_1+n_2)\bda'_1+\bda_2(m_1+m_2)+\bda_2\bda'_1}
Q_{\bda+\bda'+(3m_2+3m_1-1,3n_2+3n_1-1)}P_{\bda}P_{\bda'}\\
&-\sum\Sb \bda\in\bk_1\\ \bdb\in\bk_{-1}\endSb
q^{(n_1+n_2)\bda_1+\bdb_2(m_1+m_2)+\bdb_2\bda_1}
Q_{\bda+\bdb+(3m_2+3m_1-1,3n_2+3n_1-1)}P_{\bda}P_{\bdb})\\
=&q^{n_1m_2} \eu_{21}(m_1+m_2,n_1+n_2).
\endalign$$

$$\align  & [\eu_{23}(m_1,n_1),\eu_{32}(m_2,n_2)]\\
=&[-\mu q^{-m_1n_1}P_{-(3m_1+1,3n_1+1)}-\sum\Sb \bda\in\bk_1\\
\bdb\in\bk_{-1}\endSb q^{n_1\bdb_1+\bda_2m_1+\bda_2\bdb_1}
Q_{(3m_1+1,3n_1+1)+\bda+\bdb} P_{\bda}P_{\bdb}\\
&-\sum_{\bdb,\bdb'\in\bk_{-1}}q^{n_1\bdb'_1+\bdb_2m_1+\bdb_2\bdb'_1}Q_{(3m_1+1,3n_1+1)+\bdb+\bdb'}P_{\bdb}P_{\bdb'},
Q_{(3m_2-1,3n_2-1)}]\\
=&-\mu
q^{-m_1n_1}\delta_{(-m_1,-n_1),(m_2,n_2)}\\
&-\sum_{\bda\in\bk_1}q^{n_1m_2+\bda_2m_1+\bda_2m_2}
Q_{(3m_1+1,3n_1+1)+\bda+(3m_2-1,3n_2-1)}P_{\bda}\\
&-\sum_{\bdb\in\bk_{-1}}q^{n_1m_2+\bdb_2m_1+\bdb_2m_2}
Q_{(3m_1+1,3n_1+1)+\bdb+(3m_2-1,3n_2-1)}P_{\bdb}\\
&-\sum_{\bdb'\in\bk_{-1}}q^{n_1\bdb'_1+n_2m_1+n_2\bdb'_1}
Q_{(3m_1+1,3n_1+1)+\bdb'+(3m_2-1,3n_2-1)}P_{\bdb'}\endalign$$
$$\align
=&q^{n_1m_2}(-\sum_{\bda\in\bk_1}q^{\bda_2(m_1+m_2)}Q_{(3m_1+3m_2,3n_1+3n_2)+\bda}P_{\bda}\\
&-\sum_{\bdb\in\bk_{-1}}q^{\bdb_2(m_1+m_2)}Q_{(3m_1+3m_2,3n_1+3n_2)+\bdb}P_{\bdb}
-\frac{1}{2}\mu\delta_{(m_1+m_2,n_1+n_2),(0,0)})\\
&-q^{n_2m_1}(\sum_{\bdb\in\bk_{-1}}q^{\bdb_1(n_1+n_2)}Q_{(3m_1+3m_2,3n_1+3n_2)+\bdb}P_{\bdb}+\frac{1}{2}\mu\delta_{(m_1+m_2,n_1+n_2),(0,0)})\\
=&q^{n_1m_2}\eu_{22}(m_1+m_2,n_1+n_2)-q^{n_2m_1}\eu_{33}(m_1+m_2,n_1+n_2).
\endalign$$

$$\align  & [\eu_{23}(m_1,n_1),\eu_{33}(m_2,n_2)]\\
=&[-\mu q^{-m_1n_1}P_{-(3m_1+1,3n_1+1)}-\sum\Sb \bda\in\bk_1\\
\bdb\in\bk_{-1}\endSb q^{n_1\bdb_1+\bda_2m_1+\bda_2\bdb_1}
Q_{(3m_1+1,3n_1+1)+\bda+\bdb} P_{\bda}P_{\bdb}\\
&-\sum_{\bdb,\bdb'\in\bk_{-1}}q^{n_1\bdb'_1+\bdb_2m_1+\bdb_2\bdb'_1}Q_{(3m_1+1,3n_1+1)+\bdb+\bdb'}P_{\bdb}P_{\bdb'},\\
&\sum_{\bdb\in\bk_{-1}} q^{\bdb_1n_2} Q_{(3m_2,3n_2)+\bdb}
P_{\bdb}]\\
=&-q^{-m_1n_1+(-m_1-m-2)n_2} \mu
P_{-(3m_1+1,3n_1+1)-(3m_2,3n_2)}\\
&-\sum\Sb \bda\in\bk_1\\ \bdb\in\bk_{-1}\endSb
q^{n_1(m_2+\bdb_1)+\bda_2m_1+\bdb_2(m_2+\bdb_1)+\bdb_1n_2}
Q_{\bda+\bdb+(3m_1+1,3n_1+1)+(3m_2,3n_2)} P_{\bda}  P_{\bdb}\\
&-\sum_{\bdb,\bdb'\in\bk_{-1}}
q^{n_1(m_2+\bdb'_1)+\bdb_2m_1+\bdb_2(m_2+\bdb'_1)+\bdb'_1n_2}
Q_{\bdb+(3m_2,3n_2)+\bdb'+(3m_1+1,3n_1+1)} P_{\bdb} P_{\bdb'}\\
&-\sum_{\bdb,\bdb'\in\bk_{-1}}
q^{n_1\bdb'_1+(n_2+\bdb_2)m_1+(\bdb_2+n_2)\bdb'_1+\bdb_1n_2}
Q_{\bdb+(3m_2,3n_2)+\bdb'+(3m_1+1,3n_1+1)} P_{\bdb} P_{\bdb'}\\
&+\sum_{\bdb,\bdb'\in\bk_{-1}}
q^{(\bdb_1+\bdb'_1+m_1)n_2+n_1\bdb'_1+\bdb_2m_1+\bdb_2\bdb'_1}
Q_{\bdb+(3m_2,3n_2)+\bdb'+(3m_1+1,3n_1+1)} P_{\bdb}
P_{\bdb'}\endalign$$ \noindent (the last two terms cancel)
$$\align
=&q^{n_1m_2}(-q^{-(m_1+m_2)(n_1+n_2)}\mu
P_{-(3m_1+3m_2+1,3n_1+3n_2+1)}\\
&-\sum\Sb \bda\in\bk_1\\ \bdb\in\bk_{-1}\endSb
q^{(n_1+n_2)\bdb_1+\bda_2(m_1+m_2)+\bda_2\bdb_1}
Q_{\bda+\bdb+(3m_1+3m_2+1,3n_1+3n_2+1)} P_{\bda}  P_{\bdb}\\
&-\sum_{\bdb,\bdb'\in\bk_{-1}}
q^{(n_1+n_2)\bdb'_1+\bdb_2(m_1+m_2)+\bdb_2\bdb'_1}
Q_{\bdb+\bdb'+(3m_1+3m_2+1,3n_1+3n_2+1)} P_{\bdb} P_{\bdb'})\\
=&q^{n_1m_2}\eu_{23}(m_1+m_2,n_1+n_2).
\endalign$$

The following five brackets are easy.
$$\align  & [\eu_{31}(m_1,n_1),\eu_{31}(m_2,n_2)]=0\\
& [\eu_{31}(m_1,n_1),\eu_{32}(m_2,n_2)]=0\\
&[\eu_{31}(m_1,n_1),\eu_{33}(m_2,n_2)]=-q^{m_1n_2}\eu_{31}(m_1+m_2,n_1+n_2)\\
 & [\eu_{32}(m_1,n_1),\eu_{32}(m_2,n_2)]=0\\
  &
  [\eu_{32}(m_1,n_1),\eu_{33}(m_2,n_2)]=-q^{m_1n_2}\eu_{32}(m_1+m_2,n_1+n_2).
\endalign$$
$$\align  & [\eu_{33}(m_1,n_1),\eu_{33}(m_2,n_2)]\\
=&[\sum_{\bdb\in\bk_{-1}} q^{\bdb_1n_1} Q_{(3m_1,3n_1)+\bdb}
P_{\bdb}, \sum_{\bdb\in\bk_{-1}} q^{\bdb_1n_2}
Q_{(3m_2,3n_2)+\bdb} P_{\bdb}]\\
=&\sum_{\bdb\in\bk_{-1}} q^{(\bdb_1+m_2)n_1+\bdb_1n_2}
Q_{(3m_1,3n_1)+(3m_2,3n_2)+\bdb}P_{\bdb}\\
&-\sum_{\bdb\in\bk_{-1}}
q^{(\bdb_1+m_1)n_2+\bdb_1n_1}
Q_{(3m_1,3n_1)+(3m_2,3n_2)+\bdb}P_{\bdb}\\
=&q^{n_1m_2}\sum_{\bdb\in\bk_{-1}} q^{\bdb_1(n_1+n_2}
Q_{(3m_1,3n_1)+(3m_2,3n_2)+\bdb}P_{\bdb}\\
&-q^{m_1n_2}\sum_{\bdb\in\bk_{-1}} q^{\bdb_1(n_1+n_2}
Q_{(3m_1,3n_1)+(3m_2,3n_2)+\bdb}P_{\bdb}\\
=&q^{n_1m_2}\eu_{33}(m_1+m_2,n_1+n_2)-q^{n_2m_1}\eu_{33}(m_1+m_2,n_1+n_2).
\endalign$$

Next we check the brackets involving $D_1^{(\mu)}$ and
$D_2^{(\mu)}$.
$$\align &[D_1^{(\mu)}, \eu_{11}(m_1,n_1)]\\
=&[\sum_{\bda\in\bk_1} \bda_1 Q_{\bda}
P_{\bda}+\sum_{\bdb\in\bk_{-1}} \bdb_1 Q_{\bdb}
P_{\bdb},\sum_{\bda\in\bk_1} q^{\bda_1n_1}Q_{(3m_1, 3n_1)+\bda}
P_{\bda}+\frac{1}{2}\mu\delta_{(m_1,n_1),(0,0)}]\\
=&\sum_{\bda\in\bk_1}(m_1+\bda_1) q^{\bda_1n_1}Q_{(3m_1,
3n_1)+\bda} P_{\bda}-\sum_{\bda\in\bk_1}\bda_1
q^{\bda_1n_1}Q_{(3m_1, 3n_1)+\bda} P_{\bda}\\
=&m_1(\sum_{\bda\in\bk_1} q^{\bda_1n_1}Q_{(3m_1, 3n_1)+\bda}
P_{\bda}+\frac{1}{2}\mu\delta_{(m_1,n_1),(0,0)})\\
=&m_1\eu_{11}(m_1,n_1).
\endalign$$

The following two brackets are easy.
$$\align &[D_1^{(\mu)}, \eu_{12}(m_1,n_1)]=m_1\eu_{12}(m_1,n_1)\\
 &[D_1^{(\mu)}, \eu_{13}(m_1,n_1)]=m_1\eu_{13}(m_1,n_1).\endalign$$
$$\align &[D_1^{(\mu)}, \eu_{21}(m_1,n_1)]\\
=&[\sum_{\bda\in\bk_1} \bda_1 Q_{\bda}
P_{\bda}+\sum_{\bdb\in\bk_{-1}} \bdb_1 Q_{\bdb} P_{\bdb},-q^{-m_1n_1}\mu P_{-(3m_1-1,3n_1-1)}\\
&-\sum\Sb\bda,\bda'\in\bk_1 \endSb q^{n_1\bda_1'+ \bda_2 m_1
+\bda_2\bda_1'} Q_{\bda+\bda'+(3m_1-1, 3n_1-1)}
P_{\bda}P_{\bda'}\\
&-\sum\Sb\bda\in\bk_1\\\bdb\in\bk_{-1} \endSb q^{n_1\bda_1+ \bdb_2
m_1 +\bdb_2\bda_1} Q_{\bda+\bdb+(3m_1-1, 3n_1-1)}
P_{\bda}P_{\bdb}]\\
=&q^{-m_1n_1}\mu
(-m_1)P_{-(3m_1-1,3n_1-1)}\\
&-\sum\Sb\bda,\bda'\in\bk_1 \endSb
(\bda_1+\bda'_1+m_1)q^{n_1\bda_1'+ \bda_2 m_1 +\bda_2\bda_1'}
Q_{\bda+\bda'+(3m_1-1, 3n_1-1)} P_{\bda}P_{\bda'}\\
&+\sum\Sb\bda,\bda'\in\bk_1 \endSb \bda'_1 q^{n_1\bda_1'+ \bda_2
m_1 +\bda_2\bda_1'} Q_{\bda+\bda'+(3m_1-1, 3n_1-1)}
P_{\bda}P_{\bda'}\\
&+\sum\Sb\bda,\bda'\in\bk_1 \endSb \bda_1 q^{n_1\bda_1'+ \bda_2
m_1 +\bda_2\bda_1'} Q_{\bda+\bda'+(3m_1-1, 3n_1-1)}
P_{\bda}P_{\bda'}\\
&-\sum\Sb\bda\in\bk_1\\\bdb\in\bk_{-1} \endSb
(\bda_1+\bdb_1+m_1)q^{n_1\bda_1+ \bdb_2 m_1 +\bdb_2\bda_1}
Q_{\bda+\bdb+(3m_1-1, 3n_1-1)} P_{\bda}P_{\bdb}\\
&+\sum\Sb\bda\in\bk_1\\\bdb\in\bk_{-1} \endSb \bdb_1 q^{n_1\bda_1+
\bdb_2 m_1 +\bdb_2\bda_1} Q_{\bda+\bdb+(3m_1-1, 3n_1-1)}
P_{\bda}P_{\bdb}\\
&+\sum\Sb\bda\in\bk_1\\\bdb\in\bk_{-1} \endSb \bda_1 q^{n_1\bda_1+
\bdb_2 m_1 +\bdb_2\bda_1} Q_{\bda+\bdb+(3m_1-1, 3n_1-1)}
P_{\bda}P_{\bdb}\\
=&m_1(-q^{-m_1n_1}\mu P_{-(3m_1-1,3n_1-1)}\\
&-\sum\Sb\bda,\bda'\in\bk_1 \endSb q^{n_1\bda_1'+ \bda_2 m_1
+\bda_2\bda_1'} Q_{\bda+\bda'+(3m_1-1, 3n_1-1)}
P_{\bda}P_{\bda'}\\
&-\sum\Sb\bda\in\bk_1\\\bdb\in\bk_{-1} \endSb q^{n_1\bda_1+ \bdb_2
m_1 +\bdb_2\bda_1} Q_{\bda+\bdb+(3m_1-1, 3n_1-1)}
P_{\bda}P_{\bdb})\\
=&m_1\eu_{21}(m_1,n_1).
\endalign$$

$$\align &[D_1^{(\mu)}, \eu_{22}(m_1,n_1)]\\
=&[\sum_{\bda\in\bk_1} \bda_1 Q_{\bda}
P_{\bda}+\sum_{\bdb\in\bk_{-1}} \bdb_1 Q_{\bdb} P_{\bdb},
-\sum_{\bda\in\bk_1} q^{\bda_2m_1}
Q_{(3m_1, 3n_1)+\bda} P_{\bda}\\
&-\sum_{\bdb\in\bk_{-1}} q^{\bdb_2m_1} Q_{(3m_1, 3n_1)+\bdb}
P_{\bdb}-\frac{1}{2}\mu\delta_{(m_1,n_1),(0,0)}]\\
=&-\sum_{\bda\in\bk_1}(\bda_1+m_1) q^{\bda_2m_1} Q_{(3m_1,
3n_1)+\bda} P_{\bda}+\sum_{\bda\in\bk_1}\bda_1 q^{\bda_2m_1}
Q_{(3m_1, 3n_1)+\bda} P_{\bda}\\
&-\sum_{\bdb\in\bk_{-1}}(m_1+\bdb_1) q^{\bdb_2m_1} Q_{(3m_1,
3n_1)+\bdb} P_{\bdb}+\sum_{\bdb\in\bk_{-1}} \bdb_1 q^{\bdb_2m_1}
Q_{(3m_1, 3n_1)+\bdb} P_{\bdb}\\
=&m_1(-\sum_{\bda\in\bk_1} q^{\bda_2m_1} Q_{(3m_1, 3n_1)+\bda}
P_{\bda} -\sum_{\bdb\in\bk_{-1}} q^{\bdb_2m_1} Q_{(3m_1,
3n_1)+\bdb}
P_{\bdb}-\frac{1}{2}\mu\delta_{(m_1,n_1),(0,0)}) \\
=& m_1\eu_{22}(m_1,n_1). \endalign$$

$$\align &[D_1^{(\mu)}, \eu_{23}(m_1,n_1)]\\
=&[\sum_{\bda\in\bk_1} \bda_1 Q_{\bda}
P_{\bda}+\sum_{\bdb\in\bk_{-1}} \bdb_1 Q_{\bdb} P_{\bdb},-q^{-m_1n_1}\mu P_{-(3m_1+1,3n_1+1)}\\
&-\sum\Sb\bda\in\bk_1\\\bdb\in\bk_{-1} \endSb q^{n_1\bdb_1+ \bda_2
m_1 +\bda_2\bdb_1} Q_{\bda+\bdb+(3m_1+1, 3n_1+1)}
P_{\bda}P_{\bdb}\\
&-\sum\Sb\bdb,\bdb'\in\bk_{-1} \endSb q^{n_1\bdb_1'+ \bdb_2 m_1
+\bdb_2\bdb_1'} Q_{\bdb+\bdb'+(3m_1+1, 3n_1+1)}
P_{\bdb}P_{\bdb'}]\\
=&q^{-m_1n_1}\mu
(-m_1)P_{-(3m_1+1,3n_1+1)}\\
&-\sum\Sb\bda\in\bk_1\\\bdb\in\bk_{-1}
\endSb (\bda_1+\bdb_1+m_1)q^{n_1\bdb_1+ \bda_2 m_1 +\bda_2\bdb_1}
Q_{\bda+\bdb+(3m_1+1, 3n_1+1)} P_{\bda}P_{\bdb}\\
&+\sum\Sb\bda\in\bk_1\\\bdb\in\bk_{-1} \endSb \bda_1 q^{n_1\bdb_1+
\bda_2 m_1 +\bda_2\bdb_1} Q_{\bda+\bdb+(3m_1+1, 3n_1+1)}
P_{\bda}P_{\bdb}\\
&+\sum\Sb\bda\in\bk_1\\\bdb\in\bk_{-1} \endSb \bdb_1 q^{n_1\bdb_1+
\bda_2 m_1 +\bda_2\bdb_1} Q_{\bda+\bdb+(3m_1+1, 3n_1+1)}
P_{\bda}P_{\bdb}\\
&-\sum\Sb\bdb,\bdb'\in\bk_{-1} \endSb(\bdb_1+\bdb'_1+m_1)
q^{n_1\bdb_1'+ \bdb_2 m_1 +\bdb_2\bdb_1'} Q_{\bdb+\bdb'+(3m_1+1,
3n_1+1)} P_{\bdb}P_{\bdb'}\endalign$$
$$\align
&+\sum\Sb\bdb,\bdb'\in\bk_{-1}
\endSb \bdb_1 q^{n_1\bdb_1'+ \bdb_2 m_1 +\bdb_2\bdb_1'}
Q_{\bdb+\bdb'+(3m_1+1, 3n_1+1)}
P_{\bdb}P_{\bdb'}\\
&+\sum\Sb\bdb,\bdb'\in\bk_{-1} \endSb \bdb'_1 q^{n_1\bdb_1'+
\bdb_2 m_1 +\bdb_2\bdb_1'} Q_{\bdb+\bdb'+(3m_1+1, 3n_1+1)}
P_{\bdb}P_{\bdb'}\\
=&m_1(-q^{-m_1n_1}\mu P_{-(3m_1+1,3n_1+1)}\\
&-\sum\Sb\bda\in\bk_1\\\bdb\in\bk_{-1} \endSb q^{n_1\bdb_1+ \bda_2
m_1 +\bda_2\bdb_1} Q_{\bda+\bdb+(3m_1+1, 3n_1+1)}
P_{\bda}P_{\bdb}\\
&-\sum\Sb\bdb,\bdb'\in\bk_{-1} \endSb q^{n_1\bdb_1'+ \bdb_2 m_1
+\bdb_2\bdb_1'} Q_{\bdb+\bdb'+(3m_1+1, 3n_1+1)}
P_{\bdb}P_{\bdb'})\\
=&m_1 \eu_{23}(m_1,n_1).
\endalign$$

The following two brackets are easy.
$$\align &[D_1^{(\mu)}, \eu_{31}(m_1,n_1)]=m_1\eu_{31}(m_1,n_1)\\
&[D_1^{(\mu)}, \eu_{32}(m_1,n_1)]=m_1\eu_{32}(m_1,n_1).
\endalign$$

$$\align &[D_1^{(\mu)}, \eu_{33}(m_1,n_1)]\\
=&[\sum_{\bda\in\bk_1} \bda_1 Q_{\bda}
P_{\bda}+\sum_{\bdb\in\bk_{-1}} \bdb_1 Q_{\bdb}
P_{\bdb},\sum_{\bdb\in\bk_{-1}} q^{\bdb_1n_1}Q_{(3m_1, 3n_1)+\bdb}
P_{\bdb}+\frac{1}{2}\mu\delta_{(m_1,n_1),(0,0)}]\\
=&\sum_{\bdb\in\bk_{-1}}(m_1+\bdb_1) q^{\bdb_1n_1}Q_{(3m_1,
3n_1)+\bdb} P_{\bdb}-\sum_{\bdb\in\bk_{-1}}\bdb_1
q^{\bdb_1n_1}Q_{(3m_1, 3n_1)+\bdb} P_{\bdb}\\
=&m_1(\sum_{\bdb\in\bk_{-1}} q^{\bdb_1n_1}Q_{(3m_1, 3n_1)+\bdb}
P_{\bdb}+\frac{1}{2}\mu\delta_{(m_1,n_1),(0,0)})\\
=&m_1 \eu_{33}(m_1,n_1)
\endalign$$

Similarly, we can get
$$[D_2^{(\mu)}, \eu_{ij}(m_1,n_1)]=n_1 \eu_{ij}(m_1,n_1)$$
for $1\leq i,j\leq 3$. Finally,
$$\align &[D_1^{(\mu)}, D_2^{(\mu)}]\\
=& [\sum_{\bda\in\bk_1} \bda_1 Q_{\bda}
P_{\bda}+\sum_{\bdb\in\bk_{-1}} \bdb_1 Q_{\bdb}
P_{\bdb},\sum_{\bda\in\bk_1} \bda_2 Q_{\bda}
P_{\bda}+\sum_{\bdb\in\bk_{-1}} \bdb_2 Q_{\bdb}
P_{\bdb}]\\
=&0\endalign$$

Hence $\pi: \widetilde{\frak{gl}_{3}(\bc_q)} \to End(V)$ is a Lie
algebra homomorphism.   \qed
\enddemo

\subhead \S 2. Hermitian form for
$\widetilde{\frak{gl}_{3}(\bc_q)}$-module\endsubhead

From now on we need to assume that  $|q|=1$.

Define $\omega:\widetilde{\frak{gl}_{3}(\bc_q)}\mapsto
\widetilde{\frak{gl}_{3}(\bc_q)}$ a $\br$-linear map as the
following:

$$\align \omega(\lambda x)&=\bar{\lambda}\omega (x),\forall
\lambda\in \bc,x\in \widetilde{\frak{gl}_{3}(\bc_q)} \tag 2.1 \\
\omega(E_{ij}(a))&=(-1)^{i+j}E_{ji}(\overline{a}), a \in \bc_q \tag 2.2\\
\omega(d_s)&=d_s, \, \omega(d_t)=d_t, \, \omega(c_s)=c_s, \,
\omega(c_t)=c_t \tag 2.3
   \endalign
$$
where $\br-$linear function $\bar{ }: \bc_q\rightarrow \bc_q$ is
defined as $\overline{\lambda
s^mt^n}=\bar{\lambda}t^{-n}s^{-m}=\bar{\lambda}q^{mn}s^{-m}t^{-n}$,
and $\bar{\lambda}$ is the complex conjugate, for any $\lambda\in
\bc$, and $m,n\in \bz$.

Following from [Lemma 3.4, GZ], we have

 \proclaim{Lemma 2.4}$\omega$ is an anti-linear anti-involution of
$\widetilde{\frak{gl}_{3}(\bc_q)}$.
\endproclaim

We simply write $\pi(E_{ij}(r)).v$ as $E_{ij}(r).v$, for any $v\in
V, r\in \bc_q$.

\medskip

 In [GZ], we define a hermitian form on the basis consisting of
monomials and then use another basis consisting of iterated module
actions on the ``highest weight vector " $1$ to determine the
condition for the form being positive definite. Here we will use
the second basis directly to define the hermitian form which is
much simpler.

\proclaim{Lemma 2.5 }
$E_{12}(\alpha_1)E_{12}(\alpha_2)...E_{12}(\alpha_k)E_{32}(\beta_1)E_{32}(\beta_2)...E_{32}(\beta_l).1$
(We shall call it in level $(k,l)$ in $W$), here $k,l\in
\bz_+\bigcup\{0\}$,
$\alpha_i=s^{m_i}t^{n_i},i=1...k,\beta_j=s^{u_j}t^{v_j}$,
$j=1...l,m_i,n_i,u_j,v_j\in \bz$ forms a basis for $V$.
\endproclaim

\demo{Proof} Since
$$f_{\bda,\bdb}=\prod_{(m,n)\in\bz^2}x_{(3m+1,3n+1)}^{\bda_{(m,n)}}.
\prod_{(m',n')\in\bz^2}x_{(3m'-1,3n'-1)}^{\bdb_{(m',n')}}$$
$\bda_{(m,n)}, \bdb_{(m',n')}\in \bz_{+}\cup \{0\}$, where only
finitely many $\bda_{(m,n)}, \bdb_{(m',n')}$ are nonzero, form a
basis for $V$.

Let
$g_{\bda}=\prod_{(m,n)\in\bz^2}x_{(3m+1,3n+1)}^{\bda_{(m,n)}}$,
and
$h_{\bdb}=\prod_{(m',n')\in\bz^2}x_{(3m'-1,3n'-1)}^{\bdb_{(m',n')}}$.
It is similar to [Lemma 4.2, GZ], $g_{\bda}$ can be written as a
linear combination of $E_{12}(\alpha_1)...E_{12}(\alpha_k).1$, for
$k\leq \sum_{(m,n)}\bda_{(m,n)}$ and $h_{\bdb}$ can be written as
a linear combination of $E_{32}(\beta_1)...E_{32}(\beta_l).1$, for
$l\leq \sum_{(m',n')}\bdb_{(m',n')}$.

Since
$E_{12}(\alpha)E_{32}(\beta).u=E_{32}(\beta)E_{12}(\alpha).u$ for
any $u\in V$, $f_{(\bda,\bdb)}$ can be written as a linear
combination of
$E_{12}(\alpha_1)...E_{12}(\alpha_k)E_{32}(\beta_1)...E_{32}(\beta_l).1$.
Hence
$$E_{12}(\alpha_1)...E_{12}(\alpha_k)E_{32}(\beta_1)...E_{32}(\beta_l).1$$
form a basis for $V$.\qed
\enddemo

Let
$$\Cal{B}=\{E_{12}(\alpha_1)...E_{12}(\alpha_k)E_{32}(\beta_1)...E_{32}(\beta_l).1 |
\text{ for all }k, l \in \bn, \alpha_i,\beta_j\in \bc_q\}\tag
2.6$$ be the basis for $V$.

 \proclaim{Lemma 2.7} For any $v\in V$,
$lev(v)=lev(E_{ii}(a).v)$, $i= 1, 2, 3$;
$lev(E_{12}(a)(v))=lev(v)+(1,0)$; $lev(E_{32}(a).v)=lev(v)+(0,1)$;
$lev(E_{21}(a).v)=lev(v)-(1,0)$ or $E_{21}(a).v=0$ if
$lev(v)-(1,0)\notin \bz_+^2$; $lev(E_{23}(a).v)=lev(v)-(0,1)$ or
$E_{23}(a).v=0$ if $lev(v)-(0,1)\notin \bz_+^2$, for any $0\neq
a\in\bc_q$.
\endproclaim

\demo{Proof}  We only check those $v$ in the basis $\Cal B$.

$$\align
&E_{22}(a)E_{12}(\alpha_1)E_{12}(\alpha_2)...E_{12}(\alpha_k)E_{32}(\beta_1)E_{32}(\beta_2)...E_{32}(\beta_l).1\\
=&E_{12}(\alpha_1)E_{22}(a)E_{12}(\alpha_2)...E_{12}(\alpha_k)E_{32}(\beta_1)E_{32}(\beta_2)...E_{32}(\beta_l).1\\
&-E_{12}(\alpha_1
a)E_{12}(\alpha_2)...E_{12}(\alpha_k)E_{32}(\beta_1)E_{32}(\beta_2)...E_{32}(\beta_l).1\\
=&
E_{12}(\alpha_1)E_{12}(\alpha_2)...E_{12}(\alpha_k)E_{32}(\beta_1)E_{32}(\beta_2)...E_{32}(\beta_l).
(\frac{1}{2}\mu)\kappa(a).1\\
&-\sum\Sb i=1\endSb \Sp k\endSp E_{12}(\alpha_1)...E_{12}(\alpha_i
a)...E_{12}(\alpha_k)E_{32}(\beta_1)E_{32}(\beta_2)...E_{32}(\beta_l).1\\
&+\sum\Sb i=1\endSb \Sp l \endSp
E_{12}(\alpha_1)E_{12}(\alpha_2)...E_{12}(\alpha_k)E_{32}(\beta_1)..E_{32}(\beta_i
a)...E_{32}(\beta_l).1,
\endalign$$
so $lev(v)=lev(E_{22}(a).v)$. It is similar for $E_{11}(a),
E_{33}(a)$.

$lev(E_{12}(a)(v))=lev(v)+(1,0)$ and
$lev(E_{32}(a).v)=lev(v)+(0,1)$ are the definition of level.

For $E_{21}(a).v$, we prove by induction on the level of $v$:

$E_{21}(a).v=0$ if $lev(v)=(0,n)$, $n\in\bz_+\cup \{0\}$:

 If
$n=0$, it is obvious that $E_{21}(a).1=0$.

 Suppose it is true for $n$, then
$$\align &E_{21}(a)E_{32}(\beta_1)E_{32}(\beta_2)...E_{32}(\beta_{n+1}).1\\
=&E_{32}(\beta_1)E_{21}(a)E_{32}(\beta_2)...E_{32}(\beta_{n+1}).1-E_{31}(\beta_1
a)E_{32}(\beta_2)...E_{32}(\beta_{n+1}).1\\
=&-E_{32}(\beta_2)...E_{32}(\beta_{n+1})E_{31}(\beta_1 a).1\\
=&0\endalign$$ by induction.

Suppose $lev(E_{21}(a).v)=lev(v)-(1,0)$ or $E_{21}(a).v=0$ is true
for the $lev(v)=(m-1, n)$, then for $v=E_{12}(b)v'$ with
$lev(v')=(m-1,n)$, and $0\neq b\in\bc_q$
$$E_{21}(a).E_{12}(b)v'=E_{12}(b)E_{21}(a)v'+E_{22}(ab).v'-E_{11}(ba).v'.$$
Since $lev(E_{21}(a)v')=(m-2,n)$,
$lev(E_{21}(a).E_{12}(b)v')=(m-1,n)$ or $E_{21}(a).E_{12}(b)v'=0$.
It is similar for $E_{23}(a)$.
 \qed
\enddemo
 We easily define a contravariant (w.r.t. $\pi, \omega$)
 hermitian form on $V$ by defining on the basis $\Cal B$.

 Assume that $\mu$ is a real number, define the conjugate bilinear form
 on the elements in $\Cal B$
 by induction on the level:

 $$(1,1)=1, (1,f)=0   \text{   if   }   lev (f)\neq (0,0)\tag 2.8$$

 Suppose for any $v\in\Cal B$, $(u, v)$ is defined for any $u$ such that $lev(u)=(k',l')$, with $k'+l'=r-1$,
 if $lev(u)=(k,l)$, with $k+l=r$, then there exists a $u'$ such that $lev(u')=(k-1,l)$,
 or $lev(u')=(k,l-1)$, and some $a\in \bc_q$, such that $u=E_{12}(a)u'$
 or
$u=E_{32}(a)u'$.
 Define $$\align &(E_{12}(a)u',v)=(u',\omega(E_{12}(a))v),\tag 2.9\\
 &(E_{32}(a)u',v)=(u',\omega(E_{32}(a))v).\tag 2.10\endalign$$

 \proclaim{Theorem 2.11}The conjugate bilinear form defined above is a
 hermitian form on $V$.
 \endproclaim
 \demo{Proof} We have to check that $(E_{ij}(a)u,
 v)=(u,\omega(E_{ij}(a))v)$,  for $1\leq i, j\leq 3, a\in\bc_q$, and $(D_i.u, v)=(u,
 \omega(D_i)v)$ for $i=1,2$:

 By the definition,$$\align &(E_{12}(a)u,v)=(u,\omega(E_{12}(a))v),\\
 &(E_{32}(a)u,v)=(u,\omega(E_{32}(a))v), \endalign$$
and so
$$\align(E_{13}(a)u,v)=&([E_{12}(1), E_{23}(a)]u, v)\\
=&(E_{12}(1)E_{23}(a)u,v)-(E_{23}(a)E_{12}(1)u, v)\\
=&(u,
\omega(E_{23}(a))\omega(E_{12}(1))v)-(u,\omega(E_{12}(1))\omega(E_{23}(a))v)\\
=&(u, -\omega([E_{23}(a),E_{12}(1)])v)
=(u,\omega(E_{13}(a))v).\endalign$$

Using induction on the $lev(u)$ to prove
$(E_{11}(a).u,v)=(u,\omega(E_{11}(a)).v)$:

For any $v\in \Cal B$,
$$(E_{11}(a)1,v)=\frac{1}{2}\mu\kappa(a)(1,v)=\frac{1}{2}\mu\kappa(a)\delta_{1,v}.$$

Since $lev (E_{11}(a).v) =lev (v)$ for any $v\in\Cal B$,
 $$(1, \omega(E_{11}(a)).v)=(1, E_{11}(\bar{a}).v)=\frac{1}{2}\mu\kappa(\bar{a})\delta_{1,v}.$$
 Hence $$(E_{11}(a)1,v)=(1, \omega(E_{11}(a)).v).$$
 Suppose $(E_{11}(a)u,v)=(u, \omega(E_{11}(a)).v)$ holds true for any
 $lev (u)=(l,k)$ with $l+k=r-1.$
 For $lev (u)=(l,k)$ with $l+k=r$, then $u=E_{32}(b).u'$,
 with $lev(u')=(l,k-1)$,
 $$\align(E_{11}(a)E_{32}(b).u',v)=&(E_{32}(b)E_{11}(a).u',v)\\
 =&(E_{11}(a).u',\omega(E_{32}(b)).v)\\
 =&(u',\omega(E_{11}(a))\omega(E_{32}(b)).v)\\
 =&(u',\omega(E_{32}(b))\omega(E_{11}(a)).v)\\
 =&(E_{32}(b)u',\omega(E_{11}(a)).v)\\
 =&(u,\omega(E_{11}(a)).v),\endalign$$
or $u=E_{12}(b).u'$
 with $lev(u')=(l-1,k)$,
 $$\align(E_{11}(a)E_{12}(b).u',v)=&(E_{12}(b)E_{11}(a).u',v)+([E_{11}(a),E_{12}(b)].u',v)\\
 =&(E_{11}(a).u',\omega(E_{12}(b)).v)+(u', \omega([E_{11}(a),E_{12}(b)]).v)\\
 =&(u',\omega(E_{11}(a))\omega(E_{12}(b)).v)-(u', [\omega(E_{11}(a)),\omega(E_{12}(b))].v)\\
 =&(u',\omega(E_{12}(b))\omega(E_{11}(a)).v)\\
 =&(E_{12}(b)u',\omega(E_{11}(a)).v)\\
 =&(u,\omega(E_{11}(a)).v).\endalign$$

 Thus $(E_{11}(a).u, v)=(u, \omega(E_{11}(a)).v)$; and
 $$\align(E_{22}(a).u,v)=&([E_{21}(a), E_{12}(1)].u, v)+(E_{11}(a)u,v)\\
 =&(E_{21}(a)E_{12}(1).u,v)-(E_{12}(1)E_{21}(a).u, v)+(E_{11}(a)u,v)\\
 =&(u, \omega(E_{12}(1))\omega(E_{21}(a)).v)-(u, \omega(E_{21}(a))\omega(E_{12}(1)).v)+(u,\omega(E_{11}(a)).v)\\
 =&(u, \omega ([E_{21}(a),E_{12}(1)].v))+(u,\omega(E_{11}(a)).v)\\
 =&(u,\omega(E_{22}(a)).v).\endalign$$
 It is similar for $(E_{33}(a).u, v)=(u, \omega(E_{33}(a)).v).$

  For $D_1, D_2$, we also prove by induction on the level of $u$:

 It is obvious that $(D_1.1,v)=0,$ for any $v\in\Cal B$, so $(D_1.1,1)=(1, D_1.1)=0$,
 and suppose $(1, D_1.v)=0$ is
 true for those $lev(v)=(k',l')$, with $k'+l'=r>0$, then $(1, D_1 E_{12}(s^mt^n).v)=(1,E_{12}
 D_1.v)+(1,m.v)=0$, and $(1, D_1 E_{32}(s^mt^n).v)=(1,E_{32}
 D_1.v)+(1,m.v)=0.$ So $(D_1.1,v)=(1, D_1.v).$

Suppose for any $v\in\Cal B$, $(D_1.u,v)=(u,D_1.v)$ is true for
all $lev(u)=(k',l')$ such that $k'+l'=r$, then
$$\align (D_1.E_{12}(s^mt^n).u,v)&=(E_{12}(s^mt^n)D_1.u,v)+(m.u,v)\\
&=(D_1.u,\omega(E_{12}(s^mt^n))v)+(u,m.v)\\
&=(u,D_1\omega(E_{12}(s^mt^n))v)+(u,m.v)\\
&=(u,\omega(E_{12}(s^mt^n))D_1.v)\\
&=(E_{12}(s^mt^n)u,D_1.v).\endalign$$  It is similar for
$(D_1.E_{32}(s^mt^n).u,v)=(E_{32}(s^mt^n)u,D_1.v)$.

Hence $(D_1.u,v)=(u,D_1.v)$, and so is $(D_2.u,v)=(u,D_2.v)$. Note
that $\omega(D_i)=D_i, i=1,2$. \qed
 \enddemo

 \subhead \S 3. Conditions for unitarity \endsubhead

In this section we will determine when the hermitian form given
last section is positive definite.

Let $i\in \Bbb N$, $\gamma=(\gamma_1,..., \gamma_s)$ be the
$s-partition$ of $i$. We denote $Par_s(i)$ be the set of all
$s-partition$ of $i$.

Let $\gamma\in Par_s(N)$, we say that $\pi_1'\times \pi_2'\in
S_N\times S_N$ is equivalent to $\pi_1\times \pi_2\in S_N\times
S_N$, where $S_N$ is the permutation group of $N$ letters, if for
all $z_1,...z_N,w_1,...,w_N\in \bc_q$,
$$\kappa(z_{\pi_1'(1)}w_{\pi_2'(1)}...z_{\pi_1'(\gamma_1)}w_{\pi_2'(\gamma_1)})
...\kappa(z_{\pi_1'(\gamma_1+...\gamma_{s-1}+1)}w_{\pi_2'(\gamma_1+...+\gamma_{s-1}+1)}...
z_{\pi_1'(N)}w_{\pi_2'(N)})$$ can be obtained from the analogous
expression for $\pi_1\times \pi_2$ only by rotating the variables.
(e.g. $\kappa(z_1w_1z_2w_2z_3w_3)=\kappa(z_3w_3z_1w_1z_2w_2)$).

The following lemma is due to [JK2].

\proclaim{Lemma 3.1} Let $z_1,z_2,..z_N,w_1,w_2,..w_N \in
\bc_q[s^{\pm 1},t^{\pm 1}]$
$$\align &\bigl(\matrix
 0 & z_1 \\
  0 & 0 \\
\endmatrix\bigr)\bigl(\matrix
 0 & z_2 \\
  0 & 0 \\
\endmatrix\bigr)...\bigl(\matrix
 0 & z_N \\
  0 & 0 \\
\endmatrix\bigr)\bigl(\matrix
 0 & 0 \\
 w_1  & 0 \\
\endmatrix\bigr)\bigl(\matrix
 0 & 0 \\
 w_2  & 0 \\
\endmatrix\bigr)...\bigl(\matrix
 0 & 0 \\
 w_N  & 0 \\
\endmatrix\bigr).1 \tag 3.2\\
=&\sum\Sb s=1\endSb\Sp N\endSp\sum\Sb \gamma\in
Par_s(N)\endSb\sum\Sb [\pi_1\times\pi_2]\in (S_N\times
S_N)(\gamma)\endSb(-1)^{\gamma_1-1}(-\mu)\kappa(z_{\pi_1(1)}w_{\pi_2(1)}...z_{\pi_1(\gamma_1)}w_{\pi_2(\gamma_1)})\\
&.(-1)^{\gamma_2-1}(-\mu)\kappa(z_{\pi_1(\gamma_1+1)}w_{\pi_2(\gamma_1+1)}...z_{\pi_1(\gamma_2)}w_{\pi_2(\gamma_2)}).\\
&...(-1)^{\gamma_s-1}(-\mu)\kappa(z_{\pi_1(\gamma_1+...\gamma_{s-1}+1)}w_{\pi_2(\gamma_1+...+\gamma_{s-1}+1)}...
z_{\pi_1(N)}w_{\pi_2(N)}).1
\endalign $$
\endproclaim

\proclaim{Lemma 3.3} Let $a_i, c_i,b_j, d_j \in \bc_q,
i=1,...,m,j=1,...,n $, and $R=(a_ic_j)_{m\times m}$,
$U=(b_id_j)_{n\times n}$, and set $\Lambda=\bigl(\matrix
  R & 0 \\
  0 & U \\
\endmatrix\bigr)_{(m+n)\times (m+n)}=(\lambda_{i,j})_{(m+n)\times (m+n)}$

 $$\align &E_{21}(a_1)...E_{21}(a_m)E_{23}(b_1)...E_{23}(b_n)E_{12}(c_1)...E_{12}(c_m)E_{32}(d_1)...E_{32}(d_n).1 \tag 3.4\\
 =&\sum\Sb s=1\endSb\Sp m+n\endSp\sum\Sb \gamma\in
Par_s(m+n)\endSb\sum\Sb [\pi_1\times\pi_2]\in (S_{m+n}\times
S_{m+n})(\gamma)\endSb(-1)^{\gamma_1-1}(-\mu)\kappa(\lambda_{\pi_1(1),\pi_2(1)}...\lambda_{\pi_1(\gamma_1)\pi_2(\gamma_1)})\\
&.(-1)^{\gamma_2-1}(-\mu)\kappa(\lambda_{\pi_1(\gamma_1+1),\pi_2(\gamma_1+1)}...\lambda_{\pi_1(\gamma_2),\pi_2(\gamma_2)}).\\
&...(-1)^{\gamma_s-1}(-\mu)\kappa(\lambda_{\pi_1(\gamma_1+...\gamma_{s-1}+1),\pi_2(\gamma_1+...+\gamma_{s-1}+1)}...
\lambda_{\pi_1(N),\pi_2(N)}).1
\endalign $$

 \endproclaim

\proclaim{Remark 3.5} It is easy to see that $\lambda_{i,j}$ in
every summand should be from different rows and different columns
of $\Lambda$. And if the summand of (3.4) contains some
$\lambda_{i,j}=0$, then this summand is $0$. Hence (3.4) in fact
is the sum of those $\lambda_{i,j}$ from $R$ and $U$.
\endproclaim

 \demo{Proof} Prove by induction on $n$:

  $n=0$, (3.4) is just (3.2).

   Assume (3.4) is true up to $n-1$,
$$\align
&E_{21}(a_1)...E_{21}(a_m)E_{23}(b_1)...E_{23}(b_n)E_{12}(c_1)...E_{12}(c_m)E_{32}(d_1)...E_{32}(d_n).1\\
=&E_{21}(a_1)...E_{21}(a_m)E_{23}(b_1)...E_{23}(b_{n-1})\\
&(E_{12}(c_1)E_{23}(b_n)-E_{13}(c_1b_n))
E_{12}(c_2)...E_{12}(c_m)E_{32}(d_1)...E_{32}(d_n).1\\
=&E_{21}(a_1)...E_{21}(a_m)E_{23}(b_1)...E_{23}(b_{n-1})E_{12}(c_1)E_{23}(b_n)E_{12}(c_2)...E_{12}(c_m)E_{32}(d_1)...E_{32}(d_n).1\\
&-E_{21}(a_1)...E_{21}(a_m)E_{23}(b_1)...E_{23}(b_{n-1})E_{12}(c_1)...E_{12}(c_m)(E_{12}(c_1b_nd_1)\\
&+E_{32}(d_1)E_{13}(c_1b_n))E_{32}(d_2)...E_{32}(d_n).1\\
=&E_{21}(a_1)...E_{21}(a_m)E_{23}(b_1)...E_{23}(b_{n-1})E_{12}(c_1)E_{23}(b_n)E_{12}(c_2)...E_{12}(c_m)E_{32}(d_1)...E_{32}(d_n).1\\
&+\sum\Sb i=1\endSb\Sp n\endSp
E_{21}(a_1)...E_{21}(a_m)E_{23}(b_1)...E_{23}(b_{n-1}).\\
&E_{12}(-c_1b_nd_i)E_{12}(c_2)...E_{12}(c_m)E_{32}(d_1)...\widehat{E_{32}(d_i)}...E_{32}(d_n).1\endalign$$
$$\align
=&E_{21}(a_1)...E_{21}(a_m)E_{23}(b_1)...E_{23}(b_{n-1})E_{12}(c_1)E_{12}(c_2)...E_{12}(c_m)E_{23}(b_n)E_{32}(d_1)...E_{32}(d_n).1\\
&+\sum\Sb i=1\endSb\Sp n\endSp\sum\Sb j=1\endSb\Sp m\endSp E_{21}(a_1)...E_{21}(a_m)E_{23}(b_1)...E_{23}(b_{n-1}).\\
&E_{12}(c_1)...E_{12}(-c_jb_nd_i)...E_{12}(c_m)E_{32}(d_1)...\widehat{E_{32}(d_i)}...E_{32}(d_n).1\\
=&\sum\Sb i=1\endSb\Sp n\endSp\sum\Sb j>i\endSb
E_{21}(a_1)...E_{21}(a_m)E_{23}(b_1)...E_{23}(b_{n-1})E_{12}(c_1)E_{12}(c_2)...E_{12}(c_m)\\
&E_{32}(d_1)...\widehat{E_{32}(d_i)}...E_{32}(d_{j-1})E_{32}(-d_ib_nd_j-d_jb_nd_i)...E_{32}(d_n).1\\
&+\sum\Sb i=1\endSb\Sp n\endSp E_{21}(a_1)...E_{21}(a_m)E_{23}(b_1)...E_{23}(b_{n-1}).\\
&E_{12}(c_1)E_{12}(c_2)...E_{12}(c_m)E_{32}(d_1)...\widehat{E_{32}(d_i)}...E_{32}(d_n)(-\mu)\kappa(b_nd_i).1\\
&+\sum\Sb i=1\endSb\Sp n\endSp\sum\Sb j=1\endSb\Sp m\endSp E_{21}(a_1)...E_{21}(a_m)E_{23}(b_1)...E_{23}(b_{n-1}).\\
&E_{12}(c_1)...E_{12}(-c_jb_nd_i)...E_{12}(c_m)E_{32}(d_1)...\widehat{E_{32}(d_i)}...E_{32}(d_n).1.\\
\endalign$$

 Using (3.4) is true for $n-1$, and expanding it we can get it is
also true for $n$.\qed
 \enddemo
 \proclaim{Lemma 3.6} The hermitian form on different level is
0.
\endproclaim
\demo{Proof} Only need to prove those elements in the basis $\Cal
B$. Let $$u=E_{12}(a_1)...E_{12}(a_m)E_{32}(b_1)...E_{32}(b_n).1,
\ \ v=E_{12}(c_1)...E_{12}(c_k)E_{32}(d_1)...E_{32}(d_l).1,$$ and
$(m,n)\neq (k,l)$.

At first we prove $(u,v)=0$ with $m=0$:

If $k=0$, we can suppose $n>l$, then
$$\align(u,v)=&(E_{32}(b_1)...E_{32}(b_n).1,
E_{32}(d_1)...E_{32}(d_l).1)\\
=&((-1)^{l}E_{23}(\bar{d_l})...E_{23}(\bar{d_1})E_{32}(b_1)...E_{32}(b_n).1,1)\endalign$$
by Lemma 2.7,
$lev(E_{23}(\bar{d_l})...E_{23}(\bar{d_1})E_{32}(b_1)...E_{32}(b_n).1)=(0,n-l)$

or
$E_{23}(\bar{d_l})...E_{23}(\bar{d_1})E_{32}(b_1)...E_{32}(b_n).1=0$,
then $(u,v)=0$.

For $k>0$,
$$\align
(u,v)=&(E_{32}(b_1)...E_{32}(b_n).1,
E_{12}(c_1)...E_{12}(c_k)E_{32}(d_1)...E_{32}(d_l).1) \\
=&
(-E_{21}(\overline{c_1})E_{32}(b_1)...E_{32}(b_n).1,E_{12}(c_2)...E_{12}(c_k)E_{32}(d_1)...E_{32}(d_l).1)\endalign$$
then from Lemma 2.7,
$-E_{21}(\overline{c_1})E_{32}(b_1)...E_{32}(b_n).1=0$, then
$(u,v)=0$.

Without loss of generality, we can assume that $m\leq k$, then
$$\align
&(u,v)\\
=&(E_{12}(a_1)...E_{12}(a_m)E_{32}(b_1)...E_{32}(b_n).1,
E_{12}(c_1)...E_{12}(c_k)E_{32}(d_1)...E_{32}(d_l).1)\\
=&(E_{32}(b_1)...E_{32}(b_n).1, (-1)^m
E_{21}(\overline{a_m})...E_{21}(\overline{a_1}).E_{12}(c_1)...E_{12}(c_k)E_{32}(d_1)...E_{32}(d_l).1).
\endalign$$

From Lemma 2.7,
$$lev(E_{21}(\overline{a_m})...E_{21}(\overline{a_1}).E_{12}(c_1)...E_{12}(c_k)E_{32}(d_1)...E_{32}(d_l).1)=(k-m,n)$$
or
$$E_{21}(\overline{a_m})...E_{21}(\overline{a_1}).E_{12}(c_1)...E_{12}(c_k)E_{32}(d_1)...E_{32}(d_l).1=0,$$
then back to the case $m=0$, we get $(u,v)=0$. \qed
\enddemo

Similarly to [Proposition 4.11, GZ], and together with Lemma 3.3,
we have

 \proclaim{Proposition 3.7} The hermitian form on the same
element $h$ in level $(m,n)$ is a polynomial of $\mu$, with the
leading term is $c(-1)^{m+n}(-\mu)^{m+n}=c\mu^{m+n}$ with some
constant $c>0$.
\endproclaim

Now we can show the following theorem.

 \proclaim{Theorem 3.8} $(\pi, V)$ is unitariazable if and only if $\mu>0$.
 \endproclaim

 \demo{Proof}
From [Theorem 4.12, GZ], the hermitian form in level $(0,n)$ and
$(m,0)$ is positive definite if and only if needs $\mu>0$.

Define
$$T_{a,b}(s^{m_1}t^{n_1}s^{m_2}t^{n_2}...s^{m_k}t^{n_k})=s^{m_1+a}t^{n_1+b}s^{m_2+a}t^{n_2+b}...s^{m_k+a}t^{n_k+b}$$
($a,b\in\bz$). Extend this operator to the linear operator
$\widetilde{T_{a,b}}$ on $V$ by
$$\align &\widetilde{T_{a,b}}(E_{12}(\alpha_1)E_{12}(\alpha_2)
...E_{12}(\alpha_k)E_{32}(\beta_1)E_{32}(\beta_2)...E_{32}(\beta_l).1)\\
=&E_{12}(T_{a,b}\alpha_1)E_{12}(T_{a,b}\alpha_2)
...E_{12}(T_{a,b}\alpha_k)E_{32}(T_{a,b}\beta_1)E_{32}(T_{a,b}\beta_2)...E_{32}(T_{a,b}\beta_l).1.\endalign$$

Following Lemma 3.3, $\widetilde{T_{a,b}}$ preserves the hermitian
form on $V$. Denote
$$\align
L_{l,r}(M,N)=Span &\{E_{12}(s^{m_1}t^{n_1})...E_{12}(s^{m_l}t^{n_l})E_{32}(s^{j_1}t^{k_1})...E_{32}(s^{j_r}t^{k_r}).1\\
&|m_i,n_i\geqq 0, i=1..l, j_\iota, k_\iota \geqq 0,\\
& \sum\Sb i=1\endSb\Sp l\endSp m_i+\sum\Sb \iota=1\endSb\Sp
r\endSp j_\iota\leqq M, \sum\Sb i=1\endSb\Sp r\endSp n_i+\sum\Sb
\iota=1\endSb\Sp r\endSp k_\iota\leqq N\}.
\endalign$$

 Since the hermitian form on different level is 0, we will prove
the unitarity by induction on the level.

For any $\mu >0$, the form is definite in level $(0,n)$ ([Theorem
4.12,GZ]), and suppose it is definite in level $(r,n)$, for those
$r<m$, and it is not definite in level $(m,n)$.

From Proposition 3.7, we know that the hermitian form restrict to
this level should be positive definite for $\mu$ big enough.
Assume it is not positive definite for some $\mu>0$, then there
exist $M,N$ such that the form restrict on $L_{m,n}(M,N)$ is not
positive definite. From Proposition 3.7, the form on
$L_{l,r}(M,N)$ varies smoothly with $\mu$. Then we can find a
$\mu_0$ at which the form is not positive definite, and for all
$\mu>\mu_0$, it is positive definite. And we write $(.,.)_{\mu}$
to be the hermitian form at $\mu$.

So the radical of the form is non-trivial at $\mu_0$, i.e there
exist a nonzero $\widetilde{h}\in L_{m,n}(M,N)$, such that for any
$h\in L_{m,n}(M,N)$ we have
$$(\widetilde{h},h)_{\mu_0}=0.$$

Therefore for any arbitrary element $h_{m-1,n}$ in
$L_{m-1,n}(M,N)$, and any $c\in\bc$, we have
$$(E_{21}(c).\widetilde{h},h_{m-1,n})_{\mu_0}=0.$$

Since the form is positive definite in level $(m-1,n)$, we have
$E_{21}(c).\widetilde{h}=0$, for any $c\in\bc$. Replacing
$\widetilde{h}$ by $\widetilde{T_{-a,-b}}(\widetilde{h})$ if
necessary, we can write
$$\widetilde{h}=\sum\Sb i=1\endSb \Sp m\endSp a_i(E_{12}(1))^ix_i$$
where $x_i=\sum
E_{12}(\alpha_{i+1})..E_{12}(\alpha_m)E_{32}(\beta_1)..E_{32}(\beta_n).1$
(here it is a finite sum), and $\alpha_i,\beta_j$ is the form
$s^lt^k$ and $l,k$ can not both be 0.

Let $i_0$ be the smallest one such that $a_{i_0}\neq 0$, then
$i_0\geq 1$.

Since
 $$\align
 &E_{21}(c)(E_{12}(1))^{i}E_{12}(\alpha_{i+1})..E_{12}(\alpha_m)E_{32}(\beta_1)..E_{32}(\beta_n).1\\
 =&(E_{12}(1))^{i}E_{21}(c)E_{12}(\alpha_{i+1})..E_{12}(\alpha_m)E_{32}(\beta_1)..E_{32}(\beta_n).1\\
 &+i.(E_{12}(1))^{i-1}(E_{22}(c)-E_{11}(c))E_{12}(\alpha_{i+1})..E_{12}(\alpha_m)E_{32}(\beta_1)..E_{32}(\beta_n).1\\
 &+(-2c)\frac{i.(i-1)}{2}(E_{12}(1))^{i-1}E_{12}(\alpha_{i+1})..E_{12}(\alpha_m)E_{32}(\beta_1)..E_{32}(\beta_n).1\\
=&(E_{12}(1))^{i}E_{21}(c)E_{12}(\alpha_{i+1})..E_{12}(\alpha_m)E_{32}(\beta_1)..E_{32}(\beta_n).1\\
&+i.(E_{12}(1))^{i-1}((-2c)(m-i)-n)E_{12}(\alpha_{i+1})..E_{12}(\alpha_m)E_{32}(\beta_1)..E_{32}(\beta_n).1\\
&+i.(E_{12}(1))^{i-1}E_{12}(\alpha_{i+1})..E_{12}(\alpha_m)E_{32}(\beta_1)..E_{32}(\beta_n)(E_{22}(c)-E_{11}(c)).1\\
&+(-2c)\frac{i.(i-1)}{2}(E_{12}(1))^{i-1}E_{12}(\alpha_{i+1})..E_{12}(\alpha_m)E_{32}(\beta_1)..E_{32}(\beta_n).1\\
=&(E_{12}(1))^{i}E_{21}(c)E_{12}(\alpha_{i+1})..E_{12}(\alpha_m)E_{32}(\beta_1)..E_{32}(\beta_n).1\\
&+[ic(-\mu_0)+i((-2c)(m-i)-n)+(-2c)\frac{i.(i-1)}{2}].\\
&(E_{12}(1))^{i-1}E_{12}(\alpha_{i+1})..E_{12}(\alpha_m)E_{32}(\beta_1)..E_{32}(\beta_n).1,
 \endalign$$
we have $$E_{21}(c)\widetilde{h}=\gamma
a_{i_0}(E_{12}(1))^{i_0-1}x_{i_0}+R$$ where $R$ contains those
with power of $E_{12}(1)$ greater than $i_0-1$, and
$$\gamma=i_0c(-\mu_0)+i_0((-2c)(m-i_0)-n)+(-2c)\frac{i_0.(i_0-1)}{2}=ci_0(-\mu_0-(m-i_0)-(m-1)).$$

Since $m\geq i_0$, $i_0\geq 1$, and $\mu_0\geq 0$, $\gamma\neq 0$,
contradict with $E_{21}(c)\widetilde{h}=0$.

So for any $\mu>0$, the hermitian form is positive definite.
\qed
\enddemo

\bigskip

\heading{\bf Acknowledgments}\endheading

\medskip

I am grateful to my supervisors Professors Nantel Bergeron and Yun
Gao for their encouragement and support during the preparation of
this paper, especially to Professor Yun Gao for drawing my
attention to this subject.

\enddocument